\documentclass[leqno, final]{siamart251216}

\usepackage{graphicx}
\usepackage{subfigure}
\usepackage{epstopdf}
\usepackage{amsmath, amssymb}
\usepackage{color, xcolor}
\usepackage{algorithm}
\usepackage{multicol}
\usepackage{appendix}
\usepackage[normalem]{ulem}
\usepackage{cancel}
\usepackage{stmaryrd}
\usepackage{epstopdf}
\usepackage{microtype}
\usepackage{mathtools}
\usepackage[T1]{fontenc}
\usepackage{lmodern}

\allowdisplaybreaks[2]

\numberwithin{equation}{section}

\newtheorem{remark}{Remark}[section]

\newcommand{\norm}[1]{\left\Vert#1\right\Vert}

\newcommand{\R}{\mathbb{R}}

\newcommand{\roma}{\uppercase\expandafter{\romannumeral1}}
\newcommand{\romb}{\uppercase\expandafter{\romannumeral2}}
\newcommand{\romc}{\uppercase\expandafter{\romannumeral3}}

\usepackage{hyperref}

\usepackage{bm}

\renewcommand{\vec}[1]{\bm{#1}}
\newcommand{\mat}[1]{\mathbf{#1}}
\newcommand{\grad}{\nabla}
\newcommand{\divg}{\nabla\cdot}
\newcommand{\innerproduct}[1]{\left(#1\right)}
\newcommand{\midproduct}[1]{\left[#1\right]}
\newcommand{\jump}[1]{\left\llbracket #1\right\rrbracket}
\renewcommand{\norm}[1]{\left\lVert#1\right\rVert}
\newcommand{\enorm}[1]{\left\vert\mkern-2mu\middle\vert\mkern-2mu\middle\vert#1\middle\vert\mkern-2mu\middle\vert\mkern-2mu\right\vert}
\newcommand{\snorm}[1]{\left\lvert#1\right\rvert}
\newcommand{\coproduct}[1]{\left\langle#1\right\rangle}
\newcommand{\setprod}[1]{\left\{#1\right\}}
\renewcommand{\R}{\mathbb R}
\newcommand{\C}{\mathbb C}
\newcommand{\imagunit}{\mathbf i}
\newcommand{\dd}{\,\mathrm d}
\newcommand{\opL}{\mathcal L}
\newcommand{\opE}{\mathcal E}
\newcommand{\opT}{\mathcal T}
\newcommand{\Id}{\mathbb I}
\newcommand{\Dom}{\mathcal D}
\newcommand{\D}{\mathsf D}
\newcommand{\Rea}{\operatorname{Re}}
\newcommand{\Ima}{\operatorname{Im}}

\newsiamthm{assumption}{Assumption}

\title{\texorpdfstring{Optimal stability estimates and preasymptotic error analysis of the CIP-FEM for the Helmholtz equation with uniaxial perfectly matched layers}{Optimal stability estimates and preasymptotic error analysis of the CIP-FEM for the Helmholtz equation with uniaxial perfectly matched layers}}
\author{
	Silin Wen
	\thanks{School of Mathematics and Computational Science, Xiangtan University, Hunan 411105, China ({\tt slwen@smail.xtu.edu.cn}). This work is partially supported by the National Key R\&D Program of China grant 2024YFA1012600 and the National Natural Science Foundation of China grant 12671480. This work was supported in part by the High Performance Computing Platform of Xiangtan University.}
	\and Yu Du
	\thanks{School of Mathematics and Computational Science, Xiangtan University, Hunan 411105, China ({\tt duyu@xtu.edu.cn}). Corresponding author.}
}

\begin{document}
	\maketitle
	\vspace{-1.4in}
	\slugger{sinum}{200x}{xx}{x}{xxx--xxx}
	\vspace{1.4in}
	\setcounter{page}{1}
	\begin{abstract}
		We study the uniaxial perfectly matched layer (UPML) problem for the Helmholtz equation in two and three dimensions and establish wavenumber-explicit stability estimates at high wavenumbers. For the full-space UPML problem, we prove the optimal stability estimate with the inf-sup constant satisfying $\mu_\infty\eqsim k^{-1}$ for sufficiently large $k$ if the absorption function $\sigma\in W^{1,\infty}$ is a piecewise-power function of order at least one for which $\norm{\sigma'}_{L^\infty(\setprod{t:0<\sigma\innerproduct t\leq\delta_*})}\leq \varepsilon_*$, where $\varepsilon_*$ and $\delta_*$ depend only on the PML parameters. Using this stability estimate, we obtain the optimal inf-sup estimate $\mu_L\eqsim k^{-1}$ for the truncated UPML problem. These estimates further yield preasymptotic error estimates for a linear continuous interior penalty finite element method (CIP-FEM). Numerical simulations confirmed the theoretical predictions and demonstrated the effectiveness of CIP-FEM under high-wavenumber conditions.
	\end{abstract}
	
	
	\setcounter{page}{1}
	
	\section{Introduction}
	\label{sec:introduction}
	We consider the following Helmholtz problem in two and three dimensions ($d=2,3$)
	\begin{align}
		-\Delta u-k^2u&=f,~~\qquad\text{in }\R^d, \label{eq:helm-original}\\
		r^{\frac{d-1}{2}}\snorm{\frac{\partial u}{\partial r}-\imagunit ku}&\to 0,\qquad \text{as}~r=\snorm{\vec x}\to\infty. \label{eq:sommerfeld}
	\end{align}
	Here $k>0$ is the wavenumber and $f\in L^2(\R^d)$ is the source. For a fixed $l>0$, let $B_l=\midproduct{-l,l}^d$ denote the physical domain. We truncate the problem on a bounded computational domain by surrounding $B_l$ with a perfectly matched layer (PML) \cite{Berenger1994,ChewWeedon1994,SacksKingslandLeeLee1995,CollinoMonk1998}. Our analysis is concerned with wavenumber-explicit a priori estimates for the uniaxial PML (UPML) problem.
	
	Wavenumber-explicit stability estimates for Helmholtz problems have been developed in \cite{ChandlerWildeMonk,Spence2014,MoiolaSpence2014,DiwanMoiolaSpence,GrahamPemberySpence,MoiolaSpence2019}. For truncated radial PML problems, Li and Wu \cite{LiWu} established the optimal inf-sup estimate $\mu\eqsim k^{-1}$. For UPML problems, Chen and Xiang \cite{ChenXiang} proved that the inf-sup constant satisfies $\mu\gtrsim k^{-\frac{3}{2}}$ in two dimensions. In three dimensions, the inf-sup constant $\mu\gtrsim k^{-2}$ was proved in \cite{BaiDuWu}, which is suboptimal by a factor of $k$. In two dimensions, for piecewise constant absorption functions, Wang and Zheng \cite{WangZheng} derived the optimal inf-sup estimate $\mu\eqsim k^{-1}$ for the full-space UPML problem and its truncated version. Galkowski, Gong, Graham, Lafontaine, and Spence \cite{GalkowskiGongGrahamLafontaineSpence} showed that, for $C^\infty$ absorption functions, the inf-sup constant satisfies $\mu\gtrsim k^{-1}$. It was demonstrated in \cite{GalkowskiLafontaineSpence,AversengGalkowskiSpence} that the truncation error in radial and Cartesian PMLs decays exponentially. A systematic treatment of semiclassical analysis in the numerical analysis of Helmholtz problems is given in \cite{GalkowskiSpence2026}.
	
	The preasymptotic error analysis of finite element methods for the UPML-truncated problem hinges critically on regularity estimates for the truncated solution. For problem \eqref{eq:helm-original}--\eqref{eq:sommerfeld}, wavenumber-explicit convergence and stability theories for Galerkin and discontinuous Galerkin discretizations were established in \cite{MelenkSauter2010,MelenkSauter2011,MelenkParsaniaSauter}, and were extended to lossy Helmholtz problems by Melenk, Sauter, and Torres \cite{MelenkSauterTorres}; Lafontaine, Spence, and Wunsch \cite{LafontaineSpenceWunsch} proved quasioptimality of the $hp$-FEM for the full-space heterogeneous Helmholtz equation. Wu \cite{Wu2014} derived preasymptotic error estimates for the linear FEM and CIP-FEM; the error analysis was extended to higher-order approximations in \cite{ZhuWu,DuWu2015}. For PML-truncated problems, preasymptotic error estimates for the linear CIP-FEM and the FEM with a radial PML were given in \cite{LiWu}.
	
	In this paper, we study the UPML problem \eqref{eq:UPMLstrong} for the Helmholtz equation \eqref{eq:helm-original}. We assume that the absorption function satisfies $\sigma\in W^{1,\infty}$. Suppose that there is a fixed $\delta_*\in\innerproduct{0,\sigma_0}$ such that $\norm{\sigma'}_{L^\infty\innerproduct{\setprod{0<\sigma\le\delta_*}}}\leq\varepsilon_*$, where $\varepsilon_*>0$ depends only on the PML parameters and is independent of $k$. In particular, it is satisfied by piecewise-power absorption functions of order greater than one and by piecewise linear ones with a sufficiently small slope; in contrast to \cite{WangZheng,GalkowskiGongGrahamLafontaineSpence}, neither piecewise constancy nor $C^\infty$-smoothness is required. We prove that there exists a constant $k_0>0$ such that, for all $k\geq k_0$, the following inf-sup condition holds:
	\begin{align}
		\mu_\infty:= \inf_{0\ne u\in V}\sup_{0\ne v\in V}\frac{\snorm{a\innerproduct{u,v}}}{\enorm u\enorm v} \eqsim k^{-1}. \label{eq:main-infsup}
	\end{align}
	Here $\enorm v=\innerproduct{\norm{\grad v}_{L^2}^2+k^2\norm v_{L^2}^2}^{\frac{1}{2}}$ is the energy norm. We also establish $L^2$-norm stability for the full-space UPML problem. We propose a decomposition of the UPML operator together with a smoothing operator and a commutator, and use their properties to prove the optimal inf-sup estimate. For the truncated UPML problem, we combine the full-space stability estimate with an odd reflection argument and with decay estimates for the Green function to derive the optimal inf-sup constant $\mu_L\eqsim k^{-1}$. We further establish preasymptotic error estimates for the linear CIP-FEM and the standard FEM.
	
	The remainder of the paper is organized as follows. Section \ref{sec:model} introduces the UPML problem and proves the high-wavenumber $L^2$ a priori estimate in the full space. Section \ref{sec:resolvent} derives the optimal inf-sup estimate for the full-space problem. Section \ref{sec:truncated} establishes the corresponding estimate for the truncated problem. Section \ref{sec:cip-fem} presents a preasymptotic error analysis of the linear CIP-FEM. Section \ref{sec:numerics} presents numerical experiments confirming the theoretical estimates and illustrating the performance of the CIP-FEM at high wavenumber.
	
	We use the following notation. The relations $A\lesssim B$, $A\gtrsim B$, and $A\eqsim B$ mean $A\leq CB$, $A\geq CB$, and both inequalities, respectively, where $C>0$ is a constant independent of $k$ but possibly dependent on the fixed PML parameters. The inner product $\innerproduct{u,v}=\int_{\R^d}u\overline v\dd\vec x$ is linear in the first argument. The semiclassical spaces $H_{\hbar}^s$, their duals, and the Fourier multiplier notation are introduced in Section \ref{sec:model} and Appendix \ref{app:fourier}. We use standard notation for function spaces and norms, as in \cite{Folland,ReedSimon}.
	
	\section{The UPML problem in $\R^d$ and preliminaries}
	\label{sec:model}
	
	In this section, we introduce the UPML problem and collect the preliminaries needed for the proof of the optimal stability estimates.
	
	\subsection{The uniaxial PML problem}
	\label{subsec:absorption}
	We apply the uniaxial PML method to problem \eqref{eq:helm-original}--\eqref{eq:sommerfeld} by complex coordinate stretching outside the physical domain $B_l$ \cite{Berenger1994,ChewWeedon1994,CollinoMonk1998,SacksKingslandLeeLee1995}.
	
	Let $\sigma\in W^{1,\infty}(\R)$ be an absorption function satisfying $0\leq\sigma\leq\sigma_0$ and
	\begin{align*}
		\sigma(t) = 0,\ |t| \leq l;\quad \sigma(t) = \sigma_0,\ |t| \geq \bar l.
	\end{align*}
	Here $\bar l=l+D$. For $j=1,\ldots,d$, define
	\begin{align*}
		\sigma_j&=\sigma\innerproduct{x_j},\qquad \alpha_j=\alpha\innerproduct{x_j}=1+\innerproduct{\zeta+\imagunit}\sigma_j,\\
		\widetilde x_j&=\int_0^{x_j}\alpha\innerproduct t\dd t,\qquad j=1,\ldots,d.
	\end{align*}
	Applying this coordinate stretching to \eqref{eq:helm-original} yields the UPML equation
	\begin{align}
		-J^{-1}\divg\innerproduct{\mat A\grad u}-k^2u=f\quad\text{in }\R^d,\qquad u\in H^1(\R^d), \label{eq:UPMLstrong}
	\end{align}
	where the Jacobian determinant $J$ and the coefficient matrix $\mat A$ are given by
	\begin{align}
		J=\prod_{j=1}^d\alpha_j,\qquad \mat A=\operatorname{diag}\innerproduct{J\alpha_1^{-2},\ldots,J\alpha_d^{-2}}. \label{eq:coeff-AJ}
	\end{align}
	
	For simplicity, we write $\norm{v}_{H^j} := \norm{v}_{H^j(\R^d)}$ throughout the paper. Here and below, $(\cdot,\cdot)$ denotes the $L^2(\R^d)$ inner product. Set $V=H^1(\R^d)$ and define the sesquilinear form
	\begin{align*}
		a\innerproduct{u,v}=\innerproduct{\mat A\grad u,\grad v}-k^2\innerproduct{Ju,v},
	\end{align*}
	and the weak formulation of \eqref{eq:UPMLstrong} is finding $u\in V$ such that
	\begin{align}
		a\innerproduct{u,v}=\innerproduct{Jf,v}, \qquad v\in V. \label{eq:weak}
	\end{align}
	In what follows, the generic constant $C$ and the constants implied by $\lesssim$ may depend on the fixed parameters $d,l,D,\sigma_0,\zeta$ and on $\norm{\sigma'}_{L^\infty(\R)}$, but do not depend on $k$.
	
	In this paper, we establish the optimal inf-sup estimate under the following assumption on the UPML parameters.
	\begin{assumption}\label{ass:zeta}
		Assume that $\sigma_0>0$ and $\begin{cases} \zeta=0, & d=2, \\ \zeta\ge \frac1{\sqrt3}, & d=3. \end{cases}$
	\end{assumption}
	
	\begin{lemma}
		\label{lem:sectorial-A}
		Let $\omega_0=\arctan\frac{\sigma_0}{1+\zeta\sigma_0}$, $M=\innerproduct{\innerproduct{1+\zeta\sigma_0}^2+\sigma_0^2}^{\frac{1}{2}}$, $\vartheta_2=0$ and $\vartheta_3=\frac{\omega_0}{2}$. Under Assumption~\ref{ass:zeta}, it holds that
		\begin{align}
			\Rea\setprod{e^{-\imagunit\vartheta_d}\mat A\vec\xi\cdot\overline{\vec\xi}}\geq a_*\snorm{\vec\xi}^2,\qquad \vec\xi\in\C^d, \label{eq:Aelliptic}
		\end{align}
		where
		\begin{align}
			a_*=\begin{cases} M^{-2}, & \text{if } d=2, \\
				M^{-1}\cos\innerproduct{\frac{3\omega_0}{2}}, & \text{if } d=3. \end{cases} \label{eq:astar}
		\end{align}
		The constant $a_*$ depends only on $d$ and the fixed PML parameters.
	\end{lemma}
	\begin{proof}
		For $d=2$, \eqref{eq:coeff-AJ} gives $\mat A=\operatorname{diag}\innerproduct{\frac{\alpha_2}{\alpha_1},\frac{\alpha_1}{\alpha_2}}$. Setting $\rho_j=1+\zeta\sigma_j$, we have $\alpha_j=\rho_j+\imagunit\sigma_j$, and a direct computation gives
		\begin{align*}
			\Rea\frac{\alpha_2}{\alpha_1} =\frac{\rho_1\rho_2+\sigma_1\sigma_2}{\rho_1^2+\sigma_1^2}\geq M^{-2},\qquad \Rea\frac{\alpha_1}{\alpha_2} =\frac{\rho_1\rho_2+\sigma_1\sigma_2}{\rho_2^2+\sigma_2^2}\geq M^{-2}.
		\end{align*}
		Hence, since $\vartheta_2=0$,
		\begin{align*}
			\Rea\setprod{e^{-\imagunit\vartheta_2}\mat A\vec\xi\cdot\overline{\vec\xi}} =\Rea\frac{\alpha_2}{\alpha_1}\snorm{\xi_1}^2+\Rea\frac{\alpha_1}{\alpha_2}\snorm{\xi_2}^2 \geq M^{-2}\snorm{\vec\xi}^2,
		\end{align*}
		which proves \eqref{eq:Aelliptic} for $d=2$.
		
		For $d=3$, \eqref{eq:coeff-AJ} gives $\mat A=\operatorname{diag}\innerproduct{\frac{\alpha_m\alpha_n}{\alpha_j}}$ with $\setprod{j,m,n}=\setprod{1,2,3}$. Writing $\alpha_j=\snorm{\alpha_j}e^{\imagunit\theta_j}$ with $\theta_j=\arctan\innerproduct{\frac{\sigma_j}{1+\zeta\sigma_j}}$, we have $0\leq\theta_j\leq\omega_0$, so the phase $\phi_j=\theta_m+\theta_n-\theta_j$ of $A_{jj}$ satisfies
		\begin{align*}
			-\omega_0\leq\phi_j\leq2\omega_0,\qquad -\frac{3\omega_0}{2}\leq\phi_j-\frac{\omega_0}{2}\leq\frac{3\omega_0}{2}.
		\end{align*}
		Moreover, $\snorm{A_{jj}}=\snorm{\alpha_m}\frac{\snorm{\alpha_n}}{\snorm{\alpha_j}}\geq M^{-1}$. Since $\sigma_0>0$, $\omega_0>0$. By Assumption~\ref{ass:zeta}, $\zeta\geq\frac{1}{\sqrt3}$ for $d=3$, hence $\tan\omega_0=\frac{\sigma_0}{1+\zeta\sigma_0}<\sqrt3$, so $\omega_0<\frac{\pi}{3}$ and $\snorm{\phi_j-\frac{\omega_0}{2}}\leq\frac{3\omega_0}{2}<\frac{\pi}{2}$, which yields
		\begin{align*}
			\Rea\innerproduct{e^{-\imagunit\frac{\omega_0}{2}}A_{jj}} =\snorm{A_{jj}}\cos\innerproduct{\phi_j-\frac{\omega_0}{2}} \geq M^{-1}\cos\innerproduct{\frac{3\omega_0}{2}}>0.
		\end{align*}
		Summing over $j$ and using $\vartheta_3=\frac{\omega_0}{2}$,
		\begin{align*}
			\Rea\innerproduct{e^{-\imagunit\vartheta_3}\mat A\vec\xi\cdot\overline{\vec\xi}} =\sum_{j=1}^3\Rea\innerproduct{e^{-\imagunit\frac{\omega_0}{2}}A_{jj}}\snorm{\xi_j}^2 \geq M^{-1}\cos\innerproduct{\frac{3\omega_0}{2}}\snorm{\vec\xi}^2,
		\end{align*}
		which completes the proof.
	\end{proof}
	\begin{proposition}
		\label{prop:H2-regularity}
		Let $u\in V$ solve the UPML problem~\eqref{eq:weak}. If $f\in L^2$, then $u\in H^2$ and
		\begin{align*}
			\norm u_{H^2}\leq C_k\innerproduct{\norm f_{L^2}+\norm u_{H^1}},
		\end{align*}
		where $C_k$ may depend on $k$.
	\end{proposition}
	\begin{proof}
		Since $\sigma\in W^{1,\infty}$ and $\snorm{\alpha_j}\geq1$, \eqref{eq:coeff-AJ} shows that $\mat A,J\in W^{1,\infty}(\R^d)$. For any $\nu\in\mathbb Z^d$, set
		\begin{align*}
			Q_\nu=\prod_{j=1}^d\innerproduct{\nu_j,\nu_j+1},\qquad Q_\nu^*=\prod_{j=1}^d\innerproduct{\nu_j-\tfrac12,\nu_j+\tfrac32},
		\end{align*}
		and, by \cite[Theorem~4.16]{McLeanElliptic}, we have
		\begin{align*}
			\norm u_{H^2(Q_\nu)}^2 \leq C\innerproduct{\norm{Jf+k^2Ju}_{L^2(Q_\nu^*)}^2+\norm u_{H^1(Q_\nu^*)}^2}.
		\end{align*}
		Summing over $\nu$ and using that the cubes $Q_\nu^*$ have overlap multiplicity at most $2^d$,
		\begin{align*}
			\norm u_{H^2}^2
			&=\sum_{\nu\in\mathbb Z^d}\norm u_{H^2(Q_\nu)}^2\\
			&\leq C\sum_{\nu\in\mathbb Z^d}\innerproduct{\norm{Jf+k^2Ju}_{L^2(Q_\nu^*)}^2+\norm u_{H^1(Q_\nu^*)}^2}\\
			&\leq C2^d\innerproduct{\norm{Jf+k^2Ju}_{L^2}^2+\norm u_{H^1}^2}.
		\end{align*}
		Here $C$ depends only on $d$, $a_*$, and $\norm{\mat A}_{W^{1,\infty}(\R^d)}$.
		Hence
		\begin{align*}
			\norm u_{H^2} &\leq C\innerproduct{\norm{Jf+k^2Ju}_{L^2}+\norm u_{H^1}} \leq C_k\innerproduct{\norm f_{L^2}+\norm u_{H^1}}.
		\end{align*}
		This completes the proof.
	\end{proof}
	
	\subsection{\texorpdfstring{Decomposition of the UPML operator}{Decomposition of the UPML operator}}
	In this subsection, we rewrite the UPML problem in semiclassical form. 
	
	Let $\hbar=k^{-1}$ and let $\D_j=-\imagunit\hbar\partial_j$ denote the semiclassical derivative. For $s\in\R$, the semiclassical Sobolev norm is
	\begin{align*}
		\norm u_{H_{\hbar}^s}^2=\int_{\R^d}\innerproduct{1+\hbar^2\snorm{\vec\xi}^2}^{s}\snorm{\widehat u\innerproduct{\vec\xi}}^2\dd\vec\xi,
	\end{align*}
	where $\widehat u(\xi)=(2\pi)^{-\frac d2}\int_{\R^d}u(x)e^{-\imagunit x\cdot\xi}\dd x$ is the Fourier transform; see \cite[Section~5.2]{GalkowskiSpence2026}.
	
	For fixed $\hbar>0$, the spaces $H_{\hbar}^s$ and the classical Sobolev space $H^s$ coincide as linear spaces. For $s=1$, the norm and the seminorm satisfy
	\begin{align*}
		\norm u_{H_{\hbar}^1}^2=\norm u_{L^2}^2+\snorm u_{H_{\hbar}^1}^2,\qquad\snorm u_{H_{\hbar}^1}=\innerproduct{\sum_{j=1}^d\norm{\D_ju}_{L^2}^2}^{\frac{1}{2}}.
	\end{align*}
	We regard $H_{\hbar}^{-1}$ as the dual space of $H_{\hbar}^1$ equipped with the norm
	\begin{align}
		\norm F_{H_{\hbar}^{-1}}=\sup_{0\ne v\in H^1}\frac{\snorm{\coproduct{F,v}}}{\norm v_{H_{\hbar}^1}}. \label{eq:Hh-dual}
	\end{align}
	If $F\in L^2$, then $\coproduct{F,v}$ coincides with the $L^2$ inner product $\innerproduct{F,v}$.
	
	To expand the UPML operator and to isolate the absorption terms, we define, for $j=1,\ldots,d$,
	\begin{align*}
		R_j&=\innerproduct{1+\zeta\sigma_j}^2+\sigma_j^2, & b_j&=\frac{1}{1+\innerproduct{\zeta+\imagunit}\sigma_j}, & c_j&=\frac{\innerproduct{1+\zeta\sigma_j}^2-\sigma_j^2} {\innerproduct{\innerproduct{1+\zeta\sigma_j}^2+\sigma_j^2}^2}, \\
		d_j&=\frac{2\innerproduct{1+\zeta\sigma_j}\sigma_j} {\innerproduct{\innerproduct{1+\zeta\sigma_j}^2+\sigma_j^2}^2}, & w_j&=\frac{\sqrt{2\innerproduct{1+\zeta\sigma_j}\sigma_j}} {\innerproduct{1+\zeta\sigma_j}^2+\sigma_j^2}, & e_j&=\frac{1}{\innerproduct{1+\innerproduct{\zeta+\imagunit}\sigma_j}^2}-1.
	\end{align*}
	A direct computation gives
	\begin{align*}
		b_j^2=c_j-\imagunit d_j=1+e_j,\qquad d_j=w_j^2.
	\end{align*}
	One verifies that $b_j,e_j\in W^{1,\infty}(\R^d)$ and $w_j\in L^\infty(\R^d)\cap C^{0,\frac{1}{2}}(\R^d)$. These regularity properties, together with the coefficient estimates used below, are established in Appendix \ref{app:coefficients}.
	
	Since the stretching functions depend only on the corresponding coordinate, the UPML operator can be expanded in terms of those coordinates \cite[Section~2.1]{BaiDuWu}, that is,
	\begin{align*}
		-\hbar^2J^{-1}\partial_j\innerproduct{A_{jj}\partial_ju} &=-\hbar^2\alpha_j^{-1}\partial_j\innerproduct{\alpha_j^{-1}\partial_ju} =b_j\D_jb_j\D_ju.
	\end{align*}
	In semiclassical form, the UPML operator is defined by
	\begin{align*}
		\opL_{\hbar}=\hbar^2J^{-1}\innerproduct{-\divg\innerproduct{\mat A\grad}-k^2J} =\sum_{j=1}^db_j\D_jb_j\D_j-\Id.
	\end{align*}
	Here $\Id$ is the identity operator.
	Then \eqref{eq:UPMLstrong} reads $\opL_{\hbar}u=\hbar^2f$. Similarly, we define the Helmholtz operator
	\begin{align}
		P_{0,\hbar}=\sum_{j=1}^d\D_j^2-\Id.
	\end{align}
	
	Using the following calculations,
	\begin{align*}
		b_j\D_jb_j\D_ju &=b_j\innerproduct{\D_jb_j}\D_ju+b_j^2\D_j^2u\\
		&=b_j^2\D_j^2u-\hbar\innerproduct{1-\imagunit\zeta}\sigma_j'b_j^3\D_ju,\\
		\D_jb_j^2\D_ju &=\innerproduct{\D_jb_j^2}\D_ju+b_j^2\D_j^2u\\
		&=2b_j\innerproduct{\D_jb_j}\D_ju+b_j^2\D_j^2u\\
		&=b_j^2\D_j^2u-2\hbar\innerproduct{1-\imagunit\zeta}\sigma_j'b_j^3\D_ju,
	\end{align*}
	we obtain the decomposition of the UPML operator
	\begin{align}
		\opL_{\hbar}=P_{0,\hbar}+\opE_{\hbar}+\hbar\opT_{\hbar},\qquad \opE_{\hbar}=\sum_{j=1}^d\D_je_j\D_j,\qquad \opT_{\hbar}=\innerproduct{1-\imagunit\zeta}\sum_{j=1}^d\sigma_j'b_j^3\D_j. \label{eq:decomp}
	\end{align}
	Since $e_j\in W^{1,\infty}$ and $\sigma_j'b_j^3\in L^\infty$, the decomposition \eqref{eq:decomp} is an identity in $L^2$ for every $u\in H^2$.
	
	\section{Optimal inf-sup estimates in the full space}
	\label{sec:resolvent}
	This section is devoted to the full-space UPML problem at high wavenumbers. We establish the $L^2$ a priori estimate and source estimates, and prove the optimal inf-sup estimate $\mu_\infty\eqsim k^{-1}$.
	
	\subsection{Preliminary estimates}
	Let
	\begin{align*}
		q_\hbar\innerproduct u=\sum_{j=1}^d\norm{w_j\D_ju}_{L^2}^2
	\end{align*}
	denote the quadratic form, which measures the decay of the solution in the absorbing layer. Since $w_j\in L^\infty(\R^d)$, $q_\hbar\innerproduct u^{\frac{1}{2}}\lesssim\snorm u_{H_{\hbar}^1}$. The operator $\opE_{\hbar}$ is understood in the sense of distributions. By \eqref{eq:coeff-bounds}, $\snorm{e_j}\lesssim d_j=w_j^2$. For $u,v\in H^1$,
	\begin{align}
		\snorm{\coproduct{\opE_{\hbar}u,v}} &=\snorm{\sum_{j=1}^d\int_{\R^d}e_j\D_ju\overline{\D_jv}\dd\vec x}\notag\\
		&\leq C\sum_{j=1}^d\norm{w_j\D_ju}_{L^2}\norm{w_j\D_jv}_{L^2} \leq C_{\opE} q_\hbar\innerproduct u^{\frac{1}{2}}q_\hbar\innerproduct v^{\frac{1}{2}}. \label{eq:Eestimate}
	\end{align}
	Here $C_{\opE}$ is independent of $\hbar$. By \eqref{eq:Hh-dual},
	\begin{align}
		\norm{\opE_{\hbar}u}_{H_{\hbar}^{-1}}\leq C_{\opE} q_\hbar\innerproduct u^{\frac{1}{2}}. \label{eq:E-Hminus}
	\end{align}
	
	For $0<\delta\leq\sigma_0$, set
	\begin{align*}
		E_\delta=\setprod{t:0<\sigma\innerproduct t\leq\delta}, \qquad \Theta\innerproduct\delta=\norm{\sigma'}_{L^\infty\innerproduct{E_\delta}}.
	\end{align*}
	Since $\sigma'=0$ on $\setprod{\sigma=0}$, we have $\Theta\innerproduct{\sigma_0}=\norm{\sigma'}_{L^\infty}$. In what follows, all estimates involving \( \delta \) are stated with explicit dependence on \( \delta \).
	
	\begin{lemma}
		\label{lem:T-level}
		For any $u\in H^1$ and $0<\delta<\sigma_0$, we have
		\begin{align}
			\norm{\opT_{\hbar}u}_{L^2}\leq C_{\opT}\innerproduct{\Theta\innerproduct\delta\snorm u_{H_{\hbar}^1}+\delta^{-\frac{1}{2}}q_\hbar\innerproduct u^{\frac{1}{2}}}, \label{eq:Tlevel}
		\end{align}
		where $C_{\opT}$ is independent of $\hbar$ and $\delta$.
	\end{lemma}
	\begin{proof}
		For $0<\sigma_j\leq\delta$, we have $\snorm{\sigma_j'}\leq\Theta\innerproduct\delta$. For $\sigma_j>\delta$, Lemma \ref{lem:coefficients} gives
		\begin{align*}
			\snorm{\sigma_j'}\leq\norm{\sigma'}_{L^\infty} \leq\norm{\sigma'}_{L^\infty}\delta^{-\frac{1}{2}}\sqrt{\sigma_j} \leq\frac{\norm{\sigma'}_{L^\infty}M^2}{\sqrt2}\delta^{-\frac{1}{2}}w_j.
		\end{align*}
		Hence,
		\begin{align*}
			\snorm{\sigma_j'}\leq\Theta\innerproduct\delta+\frac{\norm{\sigma'}_{L^\infty}M^2}{\sqrt2}\delta^{-\frac{1}{2}}w_j.
		\end{align*}
		Substituting this into the expression for $\opT_{\hbar}$,
		\begin{align*}
			\norm{\opT_{\hbar}u}_{L^2} &\leq\snorm{1-\imagunit\zeta}\sum_{j=1}^d\norm{\sigma_j'b_j^3\D_ju}_{L^2}\\
			&\leq C\Theta\innerproduct\delta\sum_{j=1}^d\norm{\D_ju}_{L^2} +C\frac{\norm{\sigma'}_{L^\infty}M^2}{\sqrt2}\delta^{-\frac{1}{2}}\sum_{j=1}^d\norm{w_j\D_ju}_{L^2}\\
			&\leq C\Theta\innerproduct\delta\snorm u_{H_{\hbar}^1} +C\frac{\norm{\sigma'}_{L^\infty}M^2}{\sqrt2}\delta^{-\frac{1}{2}}q_\hbar\innerproduct u^{\frac{1}{2}}.
		\end{align*}
		This completes the proof.
	\end{proof}
	\begin{lemma}
		\label{lem:basic}
		For any $u\in H^2$, $0<\hbar\leq1$ and $0<\delta<\sigma_0$, the following estimates hold:
		\begin{align}
			\norm u_{H_{\hbar}^1}^2 &\leq C_H\innerproduct{\norm u_{L^2}^2+\norm{\opL_{\hbar}u}_{L^2}\norm u_{L^2}}, \label{eq:Hbasic}\\
			q_\hbar\innerproduct u &\leq C_q\innerproduct{\norm{\opL_{\hbar}u}_{L^2}\norm u_{L^2} +\innerproduct{\hbar\Theta\innerproduct\delta+\hbar^2\delta^{-1}}\norm u_{L^2}^2}, \label{eq:qlevel}
		\end{align}
		where the constants $C_H$ and $C_q$ are independent of $\hbar$ and $\delta$.
	\end{lemma}
	\begin{proof}
		For $u\in H^2$, by integration by parts, \eqref{eq:decomp} gives
		\begin{align}
			\innerproduct{\opL_{\hbar}u,u} =\sum_{j=1}^d\int_{\R^d}c_j\snorm{\D_ju}^2\dd\vec x-\norm u_{L^2}^2 -\imagunit q_\hbar\innerproduct u+\hbar\innerproduct{\opT_{\hbar}u,u}. \label{eq:energy}
		\end{align}
		By $\snorm{b_j}\leq1$ and the Cauchy--Schwarz inequality,
		\begin{align*}
			\norm{\opT_{\hbar}u}_{L^2} &\leq\sqrt{1+\zeta^2}\sum_{j=1}^d\norm{\sigma_j'b_j^3\D_ju}_{L^2} \leq\sqrt{1+\zeta^2}\norm{\sigma'}_{L^\infty}\sum_{j=1}^d\norm{\D_ju}_{L^2}\\
			&\leq\sqrt{d\innerproduct{1+\zeta^2}}\norm{\sigma'}_{L^\infty}\snorm u_{H_{\hbar}^1}.
		\end{align*}
		Taking the real part of \eqref{eq:energy} minus $\sigma_0$ times its imaginary part, and using \eqref{eq:coeff-bounds} and the estimate above, we obtain
		\begin{align*}
			M^{-2}\snorm u_{H_{\hbar}^1}^2 &\leq\sum_{j=1}^d\int_{\R^d}\innerproduct{c_j+\sigma_0d_j}\snorm{\D_ju}^2\dd\vec x\\
			&=\norm u_{L^2}^2+\Rea\innerproduct{\opL_{\hbar}u,u}-\sigma_0\Ima\innerproduct{\opL_{\hbar}u,u} -\hbar\Rea\innerproduct{\opT_{\hbar}u,u}+\hbar\sigma_0\Ima\innerproduct{\opT_{\hbar}u,u}\\
			&\leq\norm u_{L^2}^2+\innerproduct{1+\sigma_0}\norm{\opL_{\hbar}u}_{L^2}\norm u_{L^2} +C\hbar\snorm u_{H_{\hbar}^1}\norm u_{L^2}.
		\end{align*}
		By Young's inequality
		\begin{align*}
			C\hbar\snorm u_{H_{\hbar}^1}\norm u_{L^2}\leq\frac{1}{2M^2}\snorm u_{H_{\hbar}^1}^2+C\hbar^2\norm u_{L^2}^2,
		\end{align*}
		we obtain \eqref{eq:Hbasic}.
		
		Taking the imaginary part of \eqref{eq:energy} gives
		\begin{align*}
			q_\hbar\innerproduct u =-\Ima\innerproduct{\opL_{\hbar}u,u}+\hbar\Ima\innerproduct{\opT_{\hbar}u,u} \leq\norm{\opL_{\hbar}u}_{L^2}\norm u_{L^2}+\hbar\norm{\opT_{\hbar}u}_{L^2}\norm u_{L^2},
		\end{align*}
		and, by \eqref{eq:Tlevel},
		\begin{align*}
			q_\hbar\innerproduct u &\leq\norm{\opL_{\hbar}u}_{L^2}\norm u_{L^2} +C\hbar\Theta\innerproduct\delta\snorm u_{H_{\hbar}^1}\norm u_{L^2} +C\hbar\delta^{-\frac{1}{2}}q_\hbar\innerproduct u^{\frac{1}{2}}\norm u_{L^2}\\
			&\leq\norm{\opL_{\hbar}u}_{L^2}\norm u_{L^2} +C\hbar\Theta\innerproduct\delta\snorm u_{H_{\hbar}^1}\norm u_{L^2} +\tfrac12q_\hbar\innerproduct u+C\hbar^2\delta^{-1}\norm u_{L^2}^2.
		\end{align*}
		Using \eqref{eq:Hbasic} and Young's inequality, we obtain
		\begin{align*}
			\hbar\Theta\innerproduct\delta\snorm u_{H_{\hbar}^1}\norm u_{L^2} &\leq C_H^{\frac{1}{2}}\hbar\Theta\innerproduct\delta \innerproduct{\norm u_{L^2}^2+\norm{\opL_{\hbar}u}_{L^2}\norm u_{L^2}}^{\frac{1}{2}}\norm u_{L^2}\\
			&\leq C_H^{\frac{1}{2}}\hbar\Theta\innerproduct\delta \innerproduct{\frac{3}{2}\norm u_{L^2}^2+\frac12\norm{\opL_{\hbar}u}_{L^2}\norm u_{L^2}}\\
			&\leq C\hbar\Theta\innerproduct\delta\norm u_{L^2}^2 +C\norm{\opL_{\hbar}u}_{L^2}\norm u_{L^2}.
		\end{align*}
		The last step uses $\hbar\Theta\innerproduct\delta\leq\norm{\sigma'}_{L^\infty}$. Combining the above estimates yields \eqref{eq:qlevel}.
	\end{proof}
	
	We next introduce the smoothing operator $S_{\hbar}$ and establish a mass identity that will be used to derive the $L^2$ a priori estimate.
	
	Define
	\begin{align*}
		S_{\hbar}=\innerproduct{1-\hbar^2\Delta}^{-1},\qquad m_S\innerproduct{\vec\xi}=\frac1{1+\hbar^2\snorm{\vec\xi}^2},\qquad p_0\innerproduct{\vec\xi}=\hbar^2\snorm{\vec\xi}^2-1.
	\end{align*}
	The operators $S_{\hbar}$ and $\D_jS_{\hbar}$ are bounded and self-adjoint on $L^2$, and $P_{0,\hbar}$ is self-adjoint on $H^2$; see Appendix \ref{app:fourier} for their domains and mapping properties.
	
	\begin{lemma}\label{lm:sh}
		For any $\hbar>0$ and $F\in H_{\hbar}^{-1}$, the following estimates hold:
		\begin{align}
			\norm{S_{\hbar}}_{L^2\to L^2}\leq1,\qquad \norm{\D_jS_{\hbar}}_{L^2\to L^2}\leq\tfrac12,\qquad \norm{S_{\hbar}F}_{H_{\hbar}^1}=\norm F_{H_{\hbar}^{-1}}. \label{eq:Smaps}
		\end{align}
		Moreover,
		\begin{align}
			P_{0,\hbar}S_{\hbar}=\Id-2S_{\hbar}\quad\text{on }L^2,\qquad S_{\hbar}P_{0,\hbar}=\Id-2S_{\hbar}\quad\text{on }H^2. \label{eq:P0S-main}
		\end{align}
	\end{lemma}
	\begin{proof}
		By the Plancherel theorem, for $u\in L^2$,
		\begin{align*}
			\norm{S_{\hbar}u}_{L^2}^2 &=\int_{\R^d}\frac{\snorm{\widehat u}^2}{\innerproduct{1+\hbar^2\snorm{\vec\xi}^2}^2}\dd\vec\xi \leq\norm u_{L^2}^2,\\
			\norm{\D_jS_{\hbar}u}_{L^2}^2 &=\int_{\R^d}\frac{\hbar^2\xi_j^2}{\innerproduct{1+\hbar^2\snorm{\vec\xi}^2}^2}\snorm{\widehat u}^2\dd\vec\xi \leq\frac14\norm u_{L^2}^2.
		\end{align*}
		For $F\in H_{\hbar}^{-1}$,
		\begin{align*}
			\norm{S_{\hbar}F}_{H_{\hbar}^1}^2 &=\int_{\R^d}\innerproduct{1+\hbar^2\snorm{\vec\xi}^2}\frac{\snorm{\widehat F}^2}{\innerproduct{1+\hbar^2\snorm{\vec\xi}^2}^2}\dd\vec\xi =\int_{\R^d}\frac{\snorm{\widehat F}^2}{1+\hbar^2\snorm{\vec\xi}^2}\dd\vec\xi =\norm F_{H_{\hbar}^{-1}}^2.
		\end{align*}
		Moreover,
		\begin{align*}
			\widehat{P_{0,\hbar}S_{\hbar}u} &=\frac{\hbar^2\snorm{\vec\xi}^2-1}{1+\hbar^2\snorm{\vec\xi}^2}\widehat u =\innerproduct{1-\frac2{1+\hbar^2\snorm{\vec\xi}^2}}\widehat u =\widehat{\innerproduct{\Id-2S_{\hbar}}u},\qquad u\in L^2,\\
			\widehat{S_{\hbar}P_{0,\hbar}u} &=\frac{p_0}{1+\hbar^2\snorm{\vec\xi}^2}\widehat u =\widehat{P_{0,\hbar}S_{\hbar}u},\qquad u\in H^2,
		\end{align*}
		which proves \eqref{eq:P0S-main}.
	\end{proof}
	
	\begin{lemma}
		\label{lem:massidentity}
		For any $u\in L^2$, the following identity holds:
		\begin{align}
			\norm u_{L^2}^2=4\sum_{j=1}^d\norm{\D_jS_{\hbar}u}_{L^2}^2+\norm{\innerproduct{\Id-2S_{\hbar}}u}_{L^2}^2. \label{eq:massidentity}
		\end{align}
	\end{lemma}
	\begin{proof}
		By the Fourier multiplier representation and the Plancherel theorem,
		\begin{align*}
			4\sum_{j=1}^d\norm{\D_jS_{\hbar}u}_{L^2}^2 &=\int_{\R^d}\frac{4\hbar^2\snorm{\vec\xi}^2}{\innerproduct{1+\hbar^2\snorm{\vec\xi}^2}^2}\snorm{\widehat u\innerproduct{\vec\xi}}^2\dd\vec\xi,\\
			\norm{\innerproduct{\Id-2S_{\hbar}}u}_{L^2}^2 &=\int_{\R^d}\frac{\innerproduct{\hbar^2\snorm{\vec\xi}^2-1}^2}{\innerproduct{1+\hbar^2\snorm{\vec\xi}^2}^2}\snorm{\widehat u\innerproduct{\vec\xi}}^2\dd\vec\xi.
		\end{align*}
		Adding the two identities yields \eqref{eq:massidentity}.
	\end{proof}
	
	Set $\midproduct{X,Y}=XY-YX$. For $w_j\in L^\infty\cap C^{0,\frac{1}{2}}$, let
	\begin{align*}
		\midproduct{w_j}_{C^{0,\frac{1}{2}}}=\sup_{\vec x\ne\vec y}\frac{\snorm{w_j\innerproduct{\vec x}-w_j\innerproduct{\vec y}}}{\snorm{\vec x-\vec y}^{\frac{1}{2}}}.
	\end{align*}
	
	\begin{lemma}
		\label{lem:Sh-weight}
		For any $j=1,\ldots,d$ and $0<\hbar\leq1$, the following estimate holds:
		\begin{align}
			\norm{\midproduct{w_j,S_{\hbar}}}_{L^2\to L^2} +\sum_{m=1}^d\norm{\midproduct{w_j,\D_mS_{\hbar}}}_{L^2\to L^2} \leq C_G \hbar^{\frac{1}{2}}, \label{eq:wS-main}
		\end{align}
		where $C_G$ is independent of $\hbar$.
	\end{lemma}
	\begin{proof}
		By Lemma \ref{lm:COV}, for $u\in L^2$ and almost every $\vec x\in\R^d$,
		\begin{align*}
			\midproduct{w_j,S_{\hbar}}u\innerproduct{\vec x} &=\int_{\R^d}\innerproduct{w_j\innerproduct{\vec x}-w_j\innerproduct{\vec y}}G_\hbar\innerproduct{\vec x-\vec y}u\innerproduct{\vec y}\dd\vec y,\\
			\midproduct{w_j,\D_mS_{\hbar}}u\innerproduct{\vec x} &=\int_{\R^d}\innerproduct{w_j\innerproduct{\vec x}-w_j\innerproduct{\vec y}}\mathcal G_{m,\hbar}\innerproduct{\vec x-\vec y}u\innerproduct{\vec y}\dd\vec y.
		\end{align*}
		The definition of the Hölder seminorm and Lemma \ref{lm:c0} give
		\begin{align*}
			\snorm{\midproduct{w_j,S_{\hbar}}u\innerproduct{\vec x}} &\leq\midproduct{w_j}_{C^{0,\frac{1}{2}}}\int_{\R^d} \snorm{\vec x-\vec y}^{\frac{1}{2}}G_\hbar\innerproduct{\vec x-\vec y}\snorm{u\innerproduct{\vec y}}\dd\vec y\\
			&\leq C_0\innerproduct{\snorm{\cdot}^{\frac{1}{2}}G_\hbar*\snorm u}\innerproduct{\vec x},
		\end{align*}
		and
		\begin{align*}
			\snorm{\midproduct{w_j,\D_mS_{\hbar}}u\innerproduct{\vec x}} &\leq\midproduct{w_j}_{C^{0,\frac{1}{2}}}\int_{\R^d} \snorm{\vec x-\vec y}^{\frac{1}{2}}\snorm{\mathcal G_{m,\hbar}\innerproduct{\vec x-\vec y}}\snorm{u\innerproduct{\vec y}}\dd\vec y\\
			&\leq C_0\innerproduct{\snorm{\cdot}^{\frac{1}{2}}\snorm{\mathcal G_{m,\hbar}}*\snorm u}\innerproduct{\vec x}.
		\end{align*}
		Using Young's convolution inequality $\norm{f*g}_{L^2}\leq\norm f_{L^1}\norm g_{L^2}$ \cite[Example~1, Section~IX.4]{ReedSimon},
		\begin{align*}
			\norm{\midproduct{w_j,S_{\hbar}}u}_{L^2}& \leq C_0\norm{\snorm{\cdot}^{\frac{1}{2}}G_\hbar}_{L^1}\norm u_{L^2},\\ \norm{\midproduct{w_j,\D_mS_{\hbar}}u}_{L^2}& \leq C_0\norm{\snorm{\cdot}^{\frac{1}{2}}\snorm{\mathcal G_{m,\hbar}}}_{L^1}\norm u_{L^2}.
		\end{align*}
		We sum over $m=1,\ldots,d$. Lemma \ref{lm:c0} and the moment estimates \eqref{eq:moment-app} and \eqref{eq:derivative-moment-app} at $a=\frac{1}{2}$ then give
		\begin{align*}
			&\norm{\midproduct{w_j,S_{\hbar}}}_{L^2\to L^2} +\sum_{m=1}^d\norm{\midproduct{w_j,\D_mS_{\hbar}}}_{L^2\to L^2}\\
			&\qquad\qquad\leq C_0\int_{\R^d}\snorm{\vec z}^{\frac{1}{2}} \innerproduct{G_\hbar\innerproduct{\vec z}+\sum_{m=1}^d\snorm{\mathcal G_{m,\hbar}\innerproduct{\vec z}}}\dd\vec z \leq C_G \hbar^{\frac{1}{2}}.
		\end{align*}
	\end{proof}
	
	We next introduce the operator $B_{\hbar}$ that enters the positive commutator identity; it is a smoothed, symmetrized version of a first-order operator. Let $\eta\in C_c^\infty(\R)$ be a real cut-off function independent of $\hbar$, with $0\leq\eta\leq1$ and $\eta\innerproduct t=1$ for $\snorm t\leq\bar l$. For $j=1,\ldots,d$, we define
	\begin{align*}
		\gamma\innerproduct t=\int_0^t\eta\innerproduct s\dd s,\qquad \eta_j=\eta\innerproduct{x_j},\qquad \gamma_j=\gamma\innerproduct{x_j}.
	\end{align*}
	Then $\gamma'=\eta$, and $\gamma$ and all its derivatives are bounded. Set
	\begin{align*}
		B_{0,\hbar}&=\frac12\sum_{j=1}^d\innerproduct{\gamma_j\D_j+\D_j\gamma_j} =\sum_{j=1}^d\innerproduct{\gamma_j\D_j-\frac{\imagunit \hbar}{2}\eta_j},\\
		B_{\hbar}&=S_{\hbar}B_{0,\hbar}S_{\hbar}.
	\end{align*}
	For fixed $\hbar>0$, $B_{0,\hbar}:H^1\to L^2$ and $S_{\hbar}:L^2\to H^2$ are bounded; hence $B_{\hbar}:L^2\to H^2$ is bounded.
	
	\begin{lemma}
		\label{lem:commutators}
		The following identities hold on $H^1$ and $H^2$, respectively:
		\begin{align}
			\midproduct{\D_j,B_{\hbar}} &=-\imagunit\hbar S_{\hbar}\eta_j\D_jS_{\hbar}-\frac{\hbar^2}{2}S_{\hbar}\eta_j'S_{\hbar}, \label{eq:DB-comm}\\
			\frac{\imagunit}{\hbar}\midproduct{P_{0,\hbar},B_{\hbar}} &=S_{\hbar}\innerproduct{2\sum_{j=1}^d\D_j\eta_j\D_j-\frac{\hbar^2}{2}\sum_{j=1}^d\gamma_j'''}S_{\hbar}. \label{eq:operatorcomm-main}
		\end{align}
		Both commutators extend uniquely to bounded operators on $L^2$.
	\end{lemma}
	\begin{proof}
		Direct computation gives
		\begin{align*}
			\midproduct{\D_j,\gamma_j}=-\imagunit \hbar\eta_j,\qquad \midproduct{\D_j,\eta_j}=-\imagunit \hbar\eta_j',
		\end{align*}
		and hence
		\begin{align*}
			\midproduct{\D_j,B_{0,\hbar}} &=-\imagunit \hbar\eta_j\D_j-\frac{\hbar^2}{2}\eta_j'.
		\end{align*}
		Since $\D_j$ commutes with $S_{\hbar}$ and $B_{\hbar}=S_{\hbar}B_{0,\hbar}S_{\hbar}$, it follows that
		\begin{align*}
			\midproduct{\D_j,B_{\hbar}} &=S_{\hbar}\midproduct{\D_j,B_{0,\hbar}}S_{\hbar} =-\imagunit\hbar S_{\hbar}\eta_j\D_jS_{\hbar}-\frac{\hbar^2}{2}S_{\hbar}\eta_j'S_{\hbar}.
		\end{align*}
		Since $\gamma_j'=\eta_j$, a further computation gives
		\begin{align*}
			\midproduct{\D_j^2,\gamma_j} &=-2\imagunit \hbar\gamma_j'\D_j-\hbar^2\gamma_j'',\\
			\frac{\imagunit}{\hbar}\midproduct{\D_j^2,\frac{\gamma_j\D_j+\D_j\gamma_j}{2}} &=2\D_j\gamma_j'\D_j-\frac{\hbar^2}{2}\gamma_j'''.
		\end{align*}
		For $m\ne j$, the operators $\D_j^2$ and $\frac{\gamma_m\D_m+\D_m\gamma_m}{2}$ commute; moreover, $P_{0,\hbar}$ commutes with $S_{\hbar}$. Therefore
		\begin{align*}
			\frac{\imagunit}{\hbar}\midproduct{P_{0,\hbar},B_{\hbar}} &=\frac{\imagunit}{\hbar}S_{\hbar}\midproduct{P_{0,\hbar},B_{0,\hbar}}S_{\hbar}\\
			&=\frac{\imagunit}{\hbar}S_{\hbar}\sum_{j=1}^d \midproduct{\D_j^2,\frac{\gamma_j\D_j+\D_j\gamma_j}{2}}S_{\hbar}\\
			&=S_{\hbar}\innerproduct{2\sum_{j=1}^d\D_j\eta_j\D_j-\frac{\hbar^2}{2}\sum_{j=1}^d\gamma_j'''}S_{\hbar}.
		\end{align*}
		The identities are first verified for Schwartz functions. Their right-hand sides define bounded operators on $L^2$ by \eqref{eq:Smaps}. Since $B_{\hbar}:L^2\to H^2$ is bounded, the identities extend by density to $H^1$ and $H^2$, respectively.
	\end{proof}
	
	\begin{lemma}
		\label{lem:Bh-mapping}
		$B_{\hbar}$ is a bounded self-adjoint operator on $L^2$. For $0<\hbar\leq1$, there holds
		\begin{align}
			\norm{B_{\hbar}u}_{L^2}+\sum_{j=1}^d\norm{\D_jB_{\hbar}u}_{L^2}\leq C_{B_1}\norm u_{L^2},\qquad u\in L^2. \label{eq:Bmaps}
		\end{align}
		Furthermore, for any $u\in H^1$, one has
		\begin{align}
			q_\hbar\innerproduct{S_{\hbar}u}^{\frac{1}{2}}+q_\hbar\innerproduct{B_{\hbar}u}^{\frac{1}{2}} \leq C_{B_2}\innerproduct{q_\hbar\innerproduct u^{\frac{1}{2}}+\hbar^{\frac{1}{2}}\norm u_{H_{\hbar}^1}}, \label{eq:weightedB}
		\end{align}
		where $C_{B_1}$ and $C_{B_2}$ are independent of $\hbar$.
	\end{lemma}
	\begin{proof}
		For $0<\hbar\leq1$ and $v\in H^1$,
		\begin{align*}
			\norm{B_{0,\hbar}v}_{L^2} &=\norm{\sum_{j=1}^d\innerproduct{\gamma_j\D_jv-\frac{\imagunit \hbar}{2}\eta_jv}}_{L^2}\\
			&\leq\sum_{j=1}^d\innerproduct{\norm{\gamma_j}_{L^\infty}\norm{\D_jv}_{L^2} +\frac \hbar2\norm{\eta_j}_{L^\infty}\norm v_{L^2}} \leq C\norm v_{H_{\hbar}^1}.
		\end{align*}
		By \eqref{eq:Smaps}, for $u\in L^2$,
		\begin{align*}
			\norm{B_{0,\hbar}S_{\hbar}u}_{L^2} &\leq\sum_{j=1}^d\norm{\gamma_j}_{L^\infty}\norm{\D_jS_{\hbar}u}_{L^2} +\frac \hbar2\sum_{j=1}^d\norm{\eta_j}_{L^\infty}\norm{S_{\hbar}u}_{L^2}\\
			&\leq\innerproduct{\frac12\sum_{j=1}^d\norm{\gamma_j}_{L^\infty} +\frac \hbar2\sum_{j=1}^d\norm{\eta_j}_{L^\infty}}\norm u_{L^2} \leq C\norm u_{L^2}.
		\end{align*}
		Lemma \ref{lm:sh} therefore yields
		\begin{align}
			\norm{B_{\hbar}u}_{L^2} \leq C\norm u_{L^2}, \qquad \norm{\D_jB_{\hbar}u}_{L^2} \leq C\norm u_{L^2}.\label{eq:L2Bh}
		\end{align}
		Summing over $j$ proves \eqref{eq:Bmaps}.
		
		Integration by parts gives, for $\phi,\psi\in H^1$,
		\begin{align*}
			\innerproduct{B_{0,\hbar}\phi,\psi}=\innerproduct{\phi,B_{0,\hbar}\psi}.
		\end{align*}
		Since $S_{\hbar}$ is self-adjoint and $S_{\hbar}u,S_{\hbar}v\in H^2\subset H^1$,
		\begin{align*}
			\innerproduct{B_{\hbar}u,v} &=\innerproduct{B_{0,\hbar}S_{\hbar}u,S_{\hbar}v} =\innerproduct{S_{\hbar}u,B_{0,\hbar}S_{\hbar}v} =\innerproduct{u,B_{\hbar}v},\qquad u,v\in L^2.
		\end{align*}
		Thus, the bounded operator $B_{\hbar}$ is self-adjoint on $L^2$.
		
		To prove \eqref{eq:weightedB}, we use $\D_mS_{\hbar}=S_{\hbar}\D_m$ to write
		\begin{align*}
			B_{\hbar} &=S_{\hbar}\sum_{m=1}^d\innerproduct{\gamma_m\D_m-\frac{\imagunit \hbar}{2}\eta_m}S_{\hbar} =\sum_{m=1}^dS_{\hbar}\gamma_m\innerproduct{\D_mS_{\hbar}} -\frac{\imagunit \hbar}{2}\sum_{m=1}^dS_{\hbar}\eta_mS_{\hbar}.
		\end{align*}
		Since $w_j$ commutes with $\gamma_m$ and $\eta_m$, for $j=1,\ldots,d$ we have
		\begin{align*}
			&\midproduct{w_j,B_{\hbar}}\\ &=\sum_{m=1}^d\innerproduct{w_jS_{\hbar}\gamma_m\D_mS_{\hbar}-S_{\hbar}\gamma_m\D_mS_{\hbar}w_j} -\frac{\imagunit \hbar}{2}\sum_{m=1}^d\innerproduct{w_jS_{\hbar}\eta_mS_{\hbar}-S_{\hbar}\eta_mS_{\hbar}w_j}\\
			&=\sum_{m=1}^d\innerproduct{\midproduct{w_j,S_{\hbar}}\gamma_m\D_mS_{\hbar} +S_{\hbar}\gamma_m\midproduct{w_j,\D_mS_{\hbar}}}\\
			&\quad-\frac{\imagunit \hbar}{2}\sum_{m=1}^d\innerproduct{ \midproduct{w_j,S_{\hbar}}\eta_mS_{\hbar}+S_{\hbar}\eta_m\midproduct{w_j,S_{\hbar}}},
		\end{align*}
		and, by \eqref{eq:wS-main} and \eqref{eq:Smaps},
		\begin{align}
			\norm{\midproduct{w_j,B_{\hbar}}}_{L^2\to L^2}\leq C\hbar^{\frac{1}{2}}. \label{eq:wB-main}
		\end{align}
		The identity \eqref{eq:DB-comm} also gives
		\begin{align}
			\norm{\midproduct{\D_j,B_{\hbar}}}_{L^2\to L^2}\leq C\hbar. \label{eq:DBbound-main}
		\end{align}
		For $u\in H^1$, \eqref{eq:Smaps} and \eqref{eq:wS-main} yield
		\begin{align*}
			q_\hbar\innerproduct{S_{\hbar}u}^{\frac{1}{2}} &=\innerproduct{\sum_{j=1}^d\norm{w_j\D_jS_{\hbar}u}_{L^2}^2}^{\frac{1}{2}} =\innerproduct{\sum_{j=1}^d\norm{w_jS_{\hbar}\D_ju}_{L^2}^2}^{\frac{1}{2}}\\
			&=\innerproduct{\sum_{j=1}^d\norm{ S_{\hbar}\innerproduct{w_j\D_ju}+\midproduct{w_j,S_{\hbar}}\D_ju}_{L^2}^2}^{\frac{1}{2}}\\
			&\leq\innerproduct{\sum_{j=1}^d\innerproduct{ \norm{S_{\hbar}}_{L^2\to L^2}\norm{w_j\D_ju}_{L^2} +\norm{\midproduct{w_j,S_{\hbar}}}_{L^2\to L^2}\norm{\D_ju}_{L^2}}^2}^{\frac{1}{2}}\\
			&\leq\innerproduct{\sum_{j=1}^d\innerproduct{ \norm{w_j\D_ju}_{L^2}+C\hbar^{\frac{1}{2}}\norm{\D_ju}_{L^2}}^2}^{\frac{1}{2}}\\
			&\leq\innerproduct{\sum_{j=1}^d\norm{w_j\D_ju}_{L^2}^2}^{\frac{1}{2}} +C\hbar^{\frac{1}{2}}\innerproduct{\sum_{j=1}^d\norm{\D_ju}_{L^2}^2}^{\frac{1}{2}}\\ &=q_\hbar\innerproduct u^{\frac{1}{2}}+C\hbar^{\frac{1}{2}}\snorm u_{H_{\hbar}^1}.
		\end{align*}
		Using \eqref{eq:DB-comm} again,
		\begin{align*}
			w_j\D_jB_{\hbar}u &=w_jB_{\hbar}\D_ju+w_j\midproduct{\D_j,B_{\hbar}}u =B_{\hbar}\innerproduct{w_j\D_ju}+\midproduct{w_j,B_{\hbar}}\D_ju +w_j\midproduct{\D_j,B_{\hbar}}u.
		\end{align*}
		Combining \eqref{eq:Bmaps}, \eqref{eq:wB-main}, and \eqref{eq:DBbound-main} gives
		\begin{align*}
			q_\hbar\innerproduct{B_{\hbar}u}^{\frac{1}{2}} &=\innerproduct{\sum_{j=1}^d\norm{w_j\D_jB_{\hbar}u}_{L^2}^2}^{\frac{1}{2}}\\
			&=\innerproduct{\sum_{j=1}^d\norm{ B_{\hbar}\innerproduct{w_j\D_ju}+\midproduct{w_j,B_{\hbar}}\D_ju +w_j\midproduct{\D_j,B_{\hbar}}u}_{L^2}^2}^{\frac{1}{2}}\\
			&\leq\innerproduct{\sum_{j=1}^d\innerproduct{ C\norm{w_j\D_ju}_{L^2}+C\hbar^{\frac{1}{2}}\norm{\D_ju}_{L^2}+C\hbar\norm u_{L^2}}^2}^{\frac{1}{2}}\\
			&\leq C\innerproduct{q_\hbar\innerproduct u^{\frac{1}{2}}+\hbar^{\frac{1}{2}}\snorm u_{H_{\hbar}^1}+\hbar\norm u_{L^2}}.
		\end{align*}
		Since $\hbar\leq \hbar^{\frac{1}{2}}$, we obtain
		\begin{align*}
			q_\hbar\innerproduct{S_{\hbar}u}^{\frac{1}{2}}+q_\hbar\innerproduct{B_{\hbar}u}^{\frac{1}{2}} &\leq C\innerproduct{q_\hbar\innerproduct u^{\frac{1}{2}}+\hbar^{\frac{1}{2}}\snorm u_{H_{\hbar}^1}+\hbar\norm u_{L^2}}\\
			&\leq C_{B_2}\innerproduct{q_\hbar\innerproduct u^{\frac{1}{2}}+\hbar^{\frac{1}{2}}\norm u_{H_{\hbar}^1}}.
		\end{align*}
		This proves \eqref{eq:weightedB}.
	\end{proof}
	
	\begin{theorem}
		\label{thm:smoothing}
		For any $u\in H^2$, there holds
		\begin{align}
			\frac2\hbar\Ima\innerproduct{P_{0,\hbar}u,B_{\hbar}u} =2\sum_{j=1}^d\int_{\R^d}\eta_j\snorm{\D_jS_{\hbar}u}^2\dd\vec x -\frac{\hbar^2}{2}\sum_{j=1}^d\int_{\R^d}\gamma_j'''\snorm{S_{\hbar}u}^2\dd\vec x. \label{eq:commidentity}
		\end{align}
	\end{theorem}
	\begin{proof}
		For $u\in H^2$, $B_{\hbar}u\in H^2$. Hence,
		\begin{align*}
			\innerproduct{\midproduct{P_{0,\hbar},B_{\hbar}}u,u} &=\innerproduct{P_{0,\hbar}B_{\hbar}u,u}-\innerproduct{B_{\hbar}P_{0,\hbar}u,u}\\
			&=\innerproduct{B_{\hbar}u,P_{0,\hbar}u}-\innerproduct{P_{0,\hbar}u,B_{\hbar}u} =-2\imagunit\Ima\innerproduct{P_{0,\hbar}u,B_{\hbar}u}.
		\end{align*}
		Consequently,
		\begin{align}
			\frac{\imagunit}{\hbar}\innerproduct{\midproduct{P_{0,\hbar},B_{\hbar}}u,u} =\frac2\hbar\Ima\innerproduct{P_{0,\hbar}u,B_{\hbar}u}. \label{eq:comm-pair-left}
		\end{align}
		Integration by parts gives, for $j=1,\ldots,d$,
		\begin{align*}
			\innerproduct{\D_j\eta_j\D_jS_{\hbar}u,S_{\hbar}u} &=\innerproduct{\eta_j\D_jS_{\hbar}u,\D_jS_{\hbar}u} =\int_{\R^d}\eta_j\snorm{\D_jS_{\hbar}u}^2\dd\vec x,\\
			\innerproduct{\gamma_j'''S_{\hbar}u,S_{\hbar}u} &=\int_{\R^d}\gamma_j'''\snorm{S_{\hbar}u}^2\dd\vec x.
		\end{align*}
		Combining \eqref{eq:operatorcomm-main} and \eqref{eq:comm-pair-left}, we obtain
		\begin{align*}
			\frac2\hbar\Ima\innerproduct{P_{0,\hbar}u,B_{\hbar}u} &=\frac{\imagunit}{\hbar}\innerproduct{\midproduct{P_{0,\hbar},B_{\hbar}}u,u}\\
			&=\innerproduct{S_{\hbar}\innerproduct{2\sum_{j=1}^d\D_j\eta_j\D_j-\frac{\hbar^2}{2}\sum_{j=1}^d\gamma_j'''}S_{\hbar}u,u}\\
			&=\innerproduct{\innerproduct{2\sum_{j=1}^d\D_j\eta_j\D_j-\frac{\hbar^2}{2}\sum_{j=1}^d\gamma_j'''}S_{\hbar}u,S_{\hbar}u}\\
			&=2\sum_{j=1}^d\innerproduct{\D_j\eta_j\D_jS_{\hbar}u,S_{\hbar}u} -\frac{\hbar^2}{2}\sum_{j=1}^d\innerproduct{\gamma_j'''S_{\hbar}u,S_{\hbar}u}\\
			&=2\sum_{j=1}^d\int_{\R^d}\eta_j\snorm{\D_jS_{\hbar}u}^2\dd\vec x -\frac{\hbar^2}{2}\sum_{j=1}^d\int_{\R^d}\gamma_j'''\snorm{S_{\hbar}u}^2\dd\vec x.
		\end{align*}
		This proves \eqref{eq:commidentity}.
	\end{proof}
	
	\subsection{Wavenumber-explicit stability of the UPML problem}
	In this subsection, we combine the preceding estimates to prove the $L^2$ a priori bound of Theorem \ref{thm:L2prior-threshold}. This bound is then used to establish the optimal inf-sup estimate $\mu_\infty\eqsim k^{-1}$ for the full-space UPML problem, first for $L^2$ sources and then for general right-hand sides.
	
	We impose the condition $\Theta\innerproduct{\delta_*}<\varepsilon_*$ for some fixed $\delta_*\in\innerproduct{0,\sigma_0}$, where $\varepsilon_*>0$ is defined below.
	
	Define
	\begin{align}
		C_*=C_{\opE}C_{B_2}C_q+C_{B_1}C_{\opT}C_H^{\frac{1}{2}}+C_{\opE}C_{B_2}C_q^{\frac{1}{2}}C_H^{\frac{1}{2}}>0. \label{eq:Cstar-explicit}
	\end{align}
	Here $C_{\opT},C_H,C_q,C_{\opE},C_{B_1},C_{B_2}$ are the constants in \eqref{eq:Tlevel}, \eqref{eq:Hbasic}, \eqref{eq:qlevel}, \eqref{eq:Eestimate}, \eqref{eq:Bmaps}, and \eqref{eq:weightedB}, respectively. Thus, $C_*$ depends only on $d,l,D,\sigma_0,\zeta,\eta$, and $\norm{\sigma'}_{L^\infty}$. Set
	\begin{align}
		d_\infty=\frac{2\innerproduct{1+\zeta\sigma_0}\sigma_0} {\innerproduct{\innerproduct{1+\zeta\sigma_0}^2+\sigma_0^2}^2}\in\innerproduct{0,1}, \qquad \varepsilon_*=\innerproduct{\frac{\sqrt{1+\frac{d_\infty}{2C_*}}-1}{2}}^2>0. \label{eq:sharp-threshold}
	\end{align}
	By construction,
	\begin{align}
		C_*\innerproduct{\varepsilon_*+\varepsilon_*^{\frac{1}{2}}}=\frac{d_\infty}{8}. \label{eq:threshold-key}
	\end{align}
	
	\begin{theorem}
		\label{thm:L2prior-threshold}
		Let $\varepsilon_*$ be as in \eqref{eq:sharp-threshold}. Suppose that there exists $\delta_*\in\innerproduct{0,\sigma_0}$ such that
		\begin{align}
			\Theta\innerproduct{\delta_*}<\varepsilon_*. \label{eq:threshold-condition}
		\end{align}
		Then there exist $C>0$ and $\hbar_0\in\innerproduct{0,1}$ such that, for any $u\in H^2$ and $0<\hbar\leq \hbar_0$,
		\begin{align}
			\norm u_{L^2}\le C\hbar^{-1}\norm{\opL_{\hbar}u}_{L^2}. \label{eq:L2prior-threshold}
		\end{align}
		Here $C$ and $\hbar_0$ are independent of $\hbar$.
	\end{theorem}
	
	\begin{proof}
		Suppose that \eqref{eq:L2prior-threshold} fails. For any integer $n\geq1$, choose $C=n$ and $\hbar_0=\frac{1}{2n}$. Then there exist $0<\hbar_n\leq\frac{1}{2n}$ and $0\ne\widetilde u_n\in H^2$ such that
		\begin{align*}
			\norm{\widetilde u_n}_{L^2}>n\hbar_n^{-1}\norm{\opL_{\hbar_n}\widetilde u_n}_{L^2}.
		\end{align*}
		Set $u_n=\frac{\widetilde u_n}{\norm{\widetilde u_n}_{L^2}}$. Then
		\begin{align}
			\hbar_n\longrightarrow0,\qquad \norm{u_n}_{L^2}=1,\qquad \norm{\opL_{\hbar_n}u_n}_{L^2}<\frac{\hbar_n}{n}. \label{eq:th-contra}
		\end{align}
		By \eqref{eq:threshold-condition} and \eqref{eq:threshold-key},
		\begin{align}
			C_*\innerproduct{\Theta\innerproduct{\delta_*}+\Theta\innerproduct{\delta_*}^{\frac{1}{2}}}<\frac{d_\infty}{8}. \label{eq:th-small}
		\end{align}
		The estimates \eqref{eq:th-contra} and \eqref{eq:Hbasic} give
		\begin{align*}
			\norm{u_n}_{H_{\hbar_n}^1}^2 \leq C_H\innerproduct{1+\frac{\hbar_n}{n}}=O\innerproduct{1}.
		\end{align*}
		Keep $\delta=\delta_*$ fixed as $n\to\infty$. By \eqref{eq:qlevel} and \eqref{eq:th-contra}, $q_{\hbar_n}\innerproduct{u_n}=O(\hbar_n)$, and the first bound below holds. Moreover, \eqref{eq:Hbasic} implies $\limsup_{n\to\infty}\norm{u_n}_{H_{\hbar_n}^1}\leq C_H^{\frac{1}{2}}$. We then apply \eqref{eq:Tlevel} to obtain the second bound:
		\begin{align}
			\limsup_{n\to\infty}\frac{q_{\hbar_n}\innerproduct{u_n}}{\hbar_n} &\leq C_q\Theta\innerproduct{\delta_*}, \qquad \limsup_{n\to\infty}\norm{\opT_{\hbar_n}u_n}_{L^2} \leq C_{\opT}C_H^{\frac{1}{2}}\Theta\innerproduct{\delta_*}. \label{eq:th-fine}
		\end{align}
		The bounds \eqref{eq:weightedB} and \eqref{eq:th-fine} yield
		\begin{align*}
			q_{\hbar_n}\innerproduct{S_{\hbar_n}u_n}^{\frac{1}{2}}+q_{\hbar_n}\innerproduct{B_{\hbar_n}u_n}^{\frac{1}{2}} &\leq C_{B_2}q_{\hbar_n}\innerproduct{u_n}^{\frac{1}{2}}+C_{B_2}\hbar_n^{\frac{1}{2}}\norm{u_n}_{H_{\hbar_n}^1} =O\innerproduct{\hbar_n^{\frac{1}{2}}}.
		\end{align*}
		Hence, $q_{\hbar_n}\innerproduct{S_{\hbar_n}u_n}=o(1)$, and
		\begin{align*}
			\limsup_{n\to\infty}\hbar_n^{-\frac{1}{2}}q_{\hbar_n}\innerproduct{B_{\hbar_n}u_n}^{\frac{1}{2}} \leq C_{B_2}\innerproduct{C_q^{\frac{1}{2}}\Theta\innerproduct{\delta_*}^{\frac{1}{2}}+C_H^{\frac{1}{2}}}.
		\end{align*}
		The decomposition \eqref{eq:decomp}, \eqref{eq:Eestimate}, \eqref{eq:Bmaps}, and the preceding bounds give
		\begin{align*}
			&\frac1{\hbar_n}\snorm{\innerproduct{P_{0,\hbar_n}u_n,B_{\hbar_n}u_n}}\\ &\qquad\le\frac{\norm{\opL_{\hbar_n}u_n}_{L^2}}{\hbar_n}\norm{B_{\hbar_n}u_n}_{L^2}+\frac{C_{\opE}}{\hbar_n}q_{\hbar_n}\innerproduct{u_n}^{\frac{1}{2}}q_{\hbar_n}\innerproduct{B_{\hbar_n}u_n}^{\frac{1}{2}}\\
			&\qquad\quad+\norm{\opT_{\hbar_n}u_n}_{L^2}\norm{B_{\hbar_n}u_n}_{L^2}\\
			&\qquad\le o(1)+C_{\opE}C_{B_2}\Bigl(C_q\Theta\innerproduct{\delta_*} +C_q^{\frac{1}{2}}C_H^{\frac{1}{2}}\Theta\innerproduct{\delta_*}^{\frac{1}{2}}\Bigr) +C_{B_1}C_{\opT}C_H^{\frac{1}{2}}\Theta\innerproduct{\delta_*},
		\end{align*}
		and, by the definition of $C_*$ in \eqref{eq:Cstar-explicit},
		\begin{align*}
			\limsup_{n\to\infty}\frac1{\hbar_n}\snorm{\innerproduct{P_{0,\hbar_n}u_n,B_{\hbar_n}u_n}} \le C_*\innerproduct{\Theta\innerproduct{\delta_*}+\Theta\innerproduct{\delta_*}^{\frac{1}{2}}}.
		\end{align*}
		The positive commutator identity \eqref{eq:commidentity} gives
		\begin{align*}
			\sum_{j=1}^d\int_{\R^d}\eta_j\snorm{\D_jS_{\hbar_n}u_n}^2\dd\vec x=\frac1{\hbar_n}\Ima\innerproduct{P_{0,\hbar_n}u_n,B_{\hbar_n}u_n}+\frac{\hbar_n^2}{4}\sum_{j=1}^d\int_{\R^d}\gamma_j'''\snorm{S_{\hbar_n}u_n}^2\dd\vec x.
		\end{align*}
		The second term on the right is $o(1)$, since $\hbar_n\to0$, $\norm{S_{\hbar_n}u_n}_{L^2}\leq1$, and $\gamma_j'''\in L^\infty$. Together with $q_{\hbar_n}\innerproduct{S_{\hbar_n}u_n}=o(1)$, this yields
		\begin{align}
			\limsup_{n\to\infty}\Bigl(\sum_{j=1}^d\int_{\R^d}\eta_j\snorm{\D_jS_{\hbar_n}u_n}^2\dd\vec x +q_{\hbar_n}\innerproduct{S_{\hbar_n}u_n}\Bigr) &\le\limsup_{n\to\infty}\frac1{\hbar_n}\snorm{\innerproduct{P_{0,\hbar_n}u_n,B_{\hbar_n}u_n}}\notag\\
			&\le C_*\innerproduct{\Theta\innerproduct{\delta_*}+\Theta\innerproduct{\delta_*}^{\frac{1}{2}}} <\frac{d_\infty}{8}. \label{eq:th-eta}
		\end{align}
		The strict inequality follows from \eqref{eq:th-small}.
		
		Furthermore, \eqref{eq:decomp}, \eqref{eq:E-Hminus}, and \eqref{eq:th-fine} give the first estimate below. The second follows from \eqref{eq:P0S-main} and \eqref{eq:Smaps}:
		\begin{align*}
			\norm{P_{0,\hbar_n}u_n}_{H_{\hbar_n}^{-1}} &\le\norm{\opL_{\hbar_n}u_n}_{L^2}+Cq_{\hbar_n}\innerproduct{u_n}^{\frac{1}{2}}+\hbar_n\norm{\opT_{\hbar_n}u_n}_{L^2}=o\innerproduct{1},\\
			\norm{\innerproduct{\Id-2S_{\hbar_n}}u_n}_{L^2} &=\norm{S_{\hbar_n}P_{0,\hbar_n}u_n}_{L^2} \leq\norm{S_{\hbar_n}P_{0,\hbar_n}u_n}_{H_{\hbar_n}^1} =\norm{P_{0,\hbar_n}u_n}_{H_{\hbar_n}^{-1}}=o(1).
		\end{align*}
		For each $j$, we have $\eta_j=1$ and $d_j\geq0$ on $\setprod{\snorm{x_j}\leq\bar l}$, and $d_j=d_\infty$ and $\eta_j\geq0$ on $\setprod{\snorm{x_j}>\bar l}$. Since $0<d_\infty<1$, the bound $\eta_j+d_j\geq d_\infty$ holds in both regions. The mass identity \eqref{eq:massidentity}, the normalization $\norm{u_n}_{L^2}=1$, and \eqref{eq:th-eta} now yield
		\begin{align*}
			\frac{d_\infty}{4} &=\lim_{n\to\infty}\frac{d_\infty}{4}\innerproduct{1-\norm{\innerproduct{\Id-2S_{\hbar_n}}u_n}_{L^2}^2} =\lim_{n\to\infty}d_\infty\sum_{j=1}^d\norm{\D_jS_{\hbar_n}u_n}_{L^2}^2\\
			&=\lim_{n\to\infty}d_\infty\sum_{j=1}^d\innerproduct{\int_{\setprod{\snorm{x_j}\leq\bar l}}\snorm{\D_jS_{\hbar_n}u_n}^2\dd\vec x+\int_{\setprod{\snorm{x_j}>\bar l}}\snorm{\D_jS_{\hbar_n}u_n}^2\dd\vec x}\\
			&\leq\limsup_{n\to\infty}\sum_{j=1}^d\innerproduct{\int_{\setprod{\snorm{x_j}\leq\bar l}}\eta_j\snorm{\D_jS_{\hbar_n}u_n}^2\dd\vec x+\int_{\setprod{\snorm{x_j}>\bar l}}d_j\snorm{\D_jS_{\hbar_n}u_n}^2\dd\vec x}\\
			&\leq\limsup_{n\to\infty}\Bigl(\sum_{j=1}^d\int_{\R^d}\eta_j\snorm{\D_jS_{\hbar_n}u_n}^2\dd\vec x+q_{\hbar_n}\innerproduct{S_{\hbar_n}u_n}\Bigr) <\frac{d_\infty}{8}.
		\end{align*}
		This contradiction proves \eqref{eq:L2prior-threshold}.
	\end{proof}
	\begin{remark}
		If $\sigma(t)$ is a piecewise-power absorption function of order greater than one, the conditions of Theorem \ref{thm:L2prior-threshold} are satisfied; for a piecewise linear $\sigma(t)$, the slope must satisfy $a<\varepsilon_*$. Nevertheless, one can construct a piecewise linear absorption function of the following form:
		\begin{align*}
			\sigma\innerproduct t= \begin{cases} 0, & \snorm t\le l,\\
				a\innerproduct{\snorm t-l}, & l<\snorm t\le l+\lambda D,\\[2mm]
				\displaystyle a\lambda D+\frac{\sigma_0-a\lambda D}{\innerproduct{1-\lambda}D}\innerproduct{\snorm t-l-\lambda D}, & l+\lambda D<\snorm t<\bar l,\\[2mm]
				\sigma_0, & \snorm t\ge\bar l, \end{cases}
		\end{align*}
		with $0<\lambda<1$ and $a\lambda D<\sigma_0$.
	\end{remark}
	
	Let $V'$ be the dual space of $V$, equipped with the norm
	\begin{align*}
		\norm F_{V'}=\sup_{0\ne v\in V}\frac{\snorm{\coproduct{F,v}}}{\enorm v}.
	\end{align*}
	For $F\in V'$, the variational problem is to find $u\in V$ such that
	\begin{align}
		a\innerproduct{u,v}=\coproduct{F,v}, \qquad v\in V. \label{eq:generalweak}
	\end{align}
	The inf-sup constant associated with \eqref{eq:generalweak} is
	\begin{align}
		\mu_\infty=\inf_{0\ne u\in V}\sup_{0\ne v\in V}\frac{\snorm{a\innerproduct{u,v}}}{\enorm u\enorm v}. \label{eq:musup}
	\end{align}
	
	\begin{theorem}
		\label{thm:L2-source}
		Under the conditions of Theorem \ref{thm:L2prior-threshold}, there exist $k_0>1$ and $C>0$, both independent of $k$, such that, for any $k\geq k_0$ and $f\in L^2(\R^d)$, problem \eqref{eq:UPMLstrong} has a unique solution $u\in V$ satisfying
		\begin{align}
			k\norm u_{L^2}+\enorm u\leq C\norm f_{L^2}. \label{eq:L2-source-estimate}
		\end{align}
	\end{theorem}
	\begin{proof}
		Existence in $V$ follows from the stretched Green kernel; see \cite[Lemmas~3.2--3.3]{ChenXiang} and \cite[Section~2.1]{BaiDuWu}. Proposition \ref{prop:H2-regularity} gives $u\in H^2$. With $\hbar=k^{-1}$, Theorem \ref{thm:L2prior-threshold} and \eqref{eq:UPMLstrong} yield
		\begin{align}
			\norm u_{L^2} \leq C\hbar^{-1}\norm{\hbar^2\innerproduct{-J^{-1}\divg\innerproduct{\mat A\grad}-k^2}u}_{L^2}=Ck^{-1}\norm f_{L^2}. \label{eq:L2-source-L2}
		\end{align}
		Testing \eqref{eq:weak} with $v=u$, multiplying by $e^{-\imagunit\vartheta_d}$, and taking the real part, we obtain from \eqref{eq:Aelliptic} and \eqref{eq:L2-source-L2} that, for $k\geq k_0>1$,
		\begin{align*}
			a_*\norm{\grad u}_{L^2}^2 &\leq\snorm{\innerproduct{Jf,u}}+k^2\snorm{\innerproduct{Ju,u}}\\
			&\leq\norm J_{L^\infty}\norm f_{L^2}\norm u_{L^2}+k^2\norm J_{L^\infty}\norm u_{L^2}^2\\
			&\leq C\innerproduct{k^{-1}+1}\norm f_{L^2}^2\leq C\norm f_{L^2}^2.
		\end{align*}
		Together with \eqref{eq:L2-source-L2}, this proves \eqref{eq:L2-source-estimate}.
	\end{proof}
	\begin{theorem}
		\label{thm:high-wellposed}
		Under the conditions of Theorem \ref{thm:L2-source}, for $k\geq k_0$ and $F\in V'$, problem \eqref{eq:generalweak} has a unique solution $u\in V$. Moreover, there is a constant $C>0$ independent of $k$ such that
		\begin{align*}
			\enorm u\leq Ck\norm F_{V'},
		\end{align*}
		and $\mu_\infty\eqsim k^{-1}$.
	\end{theorem}
	\begin{proof}
		We follow the construction in \cite[Lemma~3.3]{ChenXiang} and \cite[Lemma~2.11]{BaiDuWu}. For $F\in V'$, consider the auxiliary problem
		\begin{align*}
			\innerproduct{\mat A\grad v,\grad\psi}+k^2\innerproduct{v,\psi}=\coproduct{F,\psi}, \qquad \psi\in V.
		\end{align*}
		Set
		\begin{align*}
			a_{\rm sh}=\min\setprod{a_*,\cos\vartheta_d}>0.
		\end{align*}
		By \eqref{eq:Aelliptic}, the auxiliary form multiplied by $e^{-\imagunit\vartheta_d}$ is coercive; the Lax--Milgram lemma \cite[Lemma~2.32]{McLeanElliptic} then yields a unique solution $v\in V$, and
		\begin{align*}
			a_{\rm sh}\enorm v^2\leq\Rea\setprod{e^{-\imagunit\vartheta_d}\coproduct{F,v}}\leq\norm F_{V'}\enorm v.
		\end{align*}
		Hence,
		\begin{align}
			\enorm v+k\norm v_{L^2}\leq2a_{\rm sh}^{-1}\norm F_{V'}. \label{eq:vbound}
		\end{align}
		Define $f_v=k^2J^{-1}\innerproduct{1+J}v$. By \eqref{eq:vbound},
		\begin{align}
			\norm{f_v}_{L^2}\leq k^2\norm{J^{-1}\innerproduct{1+J}}_{L^\infty}\norm v_{L^2}\leq Ck\norm F_{V'}. \label{eq:fvbound}
		\end{align}
		Let $\Phi\in V$ solve \eqref{eq:UPMLstrong} with source $f_v$. Theorem \ref{thm:L2-source} and \eqref{eq:fvbound} give
		\begin{align}
			k\norm\Phi_{L^2}+\enorm\Phi\leq C\norm{f_v}_{L^2}\leq Ck\norm F_{V'}. \label{eq:Phi-bound}
		\end{align}
		Then $u=v+\Phi$ solves \eqref{eq:generalweak}. Furthermore, \eqref{eq:vbound} and \eqref{eq:Phi-bound} imply
		\begin{align}
			\enorm u&\leq\enorm v+\enorm\Phi\leq a_{\rm sh}^{-1}\norm F_{V'}+Ck\norm F_{V'}\leq Ck\norm F_{V'}. \label{eq:dualhigh}
		\end{align}
		
		For any $u\in V$, define $F_u\in V'$ by
		\begin{align*}
			\coproduct{F_u,\psi}=a\innerproduct{u,\psi}, \qquad \psi\in V.
		\end{align*}
		Applying \eqref{eq:dualhigh} yields
		\begin{align*}
			\enorm u&\leq Ck\norm{F_u}_{V'}=Ck\sup_{0\ne\psi\in V}\frac{\snorm{a\innerproduct{u,\psi}}}{\enorm\psi},
		\end{align*}
		and therefore
		\begin{align}
			\mu_\infty\geq\frac{1}{Ck}. \label{eq:mu-lower}
		\end{align}
		The upper bound $\mu_\infty\lesssim k^{-1}$ follows from the classical localized plane-wave construction in \cite[Lemma~3.10]{ChandlerWildeMonk}; together with \eqref{eq:mu-lower} this proves $\mu_\infty\eqsim k^{-1}$.
	\end{proof}
	
	\section{The truncated UPML problem}
	\label{sec:truncated}
	We truncate the UPML problem to the bounded domain $B_L$ by imposing a homogeneous Dirichlet condition on $\partial B_L$ \cite{BramblePasciakMaxwell}. Set
	\begin{align*}
		B_L=\innerproduct{-L,L}^d, \qquad D_c=L-\bar l>0, \qquad V_L=H_0^1(B_L),
	\end{align*}
	where $D_c$ is the thickness of the outer constant absorbing layer, and $L\eqsim\bar l\eqsim l\eqsim1$. Define
	\begin{align*}
		\enorm v_{B_L}^2 &=\norm{\grad v}_{L^2(B_L)}^2 +k^2\norm v_{L^2(B_L)}^2,\qquad v\in H^1(B_L), \\
		a_L\innerproduct{u,v} &=\innerproduct{\mat A\grad u,\grad v}_{B_L} -k^2\innerproduct{Ju,v}_{B_L}, \qquad u,v\in H^1(B_L).
	\end{align*}
	For a box $B$, $a_B$ denotes the same form with integration over $B$. For $v\in H^1(B_l)$, set $\enorm v_{B_l}^2=\norm{\grad v}_{L^2(B_l)}^2+k^2\norm v_{L^2(B_l)}^2$. Let $V_L'$ be the dual space of $V_L$, equipped with the norm
	\begin{align*}
		\norm F_{V_L'}=\sup_{0\ne v\in V_L}\frac{\snorm{\coproduct{F,v}}}{\enorm v_{B_L}}.
	\end{align*}
	Set
	\begin{align}
		\mu_L =\inf_{0\ne u\in V_L}\sup_{0\ne v\in V_L} \frac{\snorm{a_L\innerproduct{u,v}}} {\enorm u_{B_L}\enorm v_{B_L}}. \label{eq:mu-truncated}
	\end{align}
	For $F\in V_L'$, the truncated problem is to find $u_L\in V_L$ such that
	\begin{align}
		a_L\innerproduct{u_L,v}=\coproduct{F,v}, \qquad v\in V_L. \label{eq:truncated-weak}
	\end{align}
	
	Since $1\leq\snorm{\alpha_j}\leq M$,
	\begin{align}
		\norm J_{L^\infty}\leq M^d, \qquad \norm{\mat A}_{L^\infty}\leq M^{d-1}. \label{eq:truncated-coefficient-bounds}
	\end{align}
	For all $u,v\in H^1(B_L)$,
	\begin{align}
		\snorm{a_L\innerproduct{u,v}}\leq\innerproduct{M^{d-1}+M^d}\enorm u_{B_L}\enorm v_{B_L}. \label{eq:truncated-continuity}
	\end{align}
	By Theorem \ref{thm:high-wellposed}, there exist $k_0>1$ and $C_\infty>0$, both independent of $k$, such that, for all $k\geq k_0$ and $u\in H^1(\R^d)$,
	\begin{align}
		\enorm u &\leq C_\infty k \sup_{0\ne v\in H^1(\R^d)} \frac{\snorm{a\innerproduct{u,v}}}{\enorm v}. \label{eq:C-infty-truncated}
	\end{align}
	Here $C_\infty$ is independent of $k$, $L$, and $D_c$.
	\begin{lemma}
		\label{lem:reflection-extension}
		Odd reflection across the coordinate planes $x_j=\pm L$ extends every $w\in V_L$ to a function $\widetilde w\in H^1(B_{L+D_c})$ satisfying
		\begin{align}
			\widetilde w|_{B_L}=w, \qquad \enorm{\widetilde w}_{B_{L+D_c}}\leq2^{\frac{d}{2}}\enorm w_{B_L}, \label{eq:reflection-extension}
		\end{align}
		and
		\begin{align}
			\sup_{0\ne v\in H_0^1(B_{L+D_c})}\frac{\snorm{a_{B_{L+D_c}}\innerproduct{\widetilde w,v}}}{\enorm v_{B_{L+D_c}}} &\leq2^{\frac{d}{2}}\sup_{0\ne v\in V_L}\frac{\snorm{a_L\innerproduct{w,v}}}{\enorm v_{B_L}}. \label{eq:reflection-pairing}
		\end{align}
	\end{lemma}
	The scalar odd reflection and its dual pairing are described in \cite[Definitions~3.1--3.2, Lemma~3.3]{GalkowskiGongGrahamLafontaineSpence}. The reflection is confined to the layer of constant absorption, where the coefficients are invariant under the corresponding reflections. Each coordinate direction contributes at most a factor of $\sqrt2$ to the norm bounds, yielding $2^{\frac{d}{2}}$ after successive reflections.
	
	We use the stretched Green function to estimate the truncation error \cite{WangZheng,BaiDuWu}. For $\vec x\ne\vec y$, define the complex distance $\rho$ and the Green kernel $G_d$ by
	\begin{align*}
		\rho\innerproduct{\widetilde{\vec x},\widetilde{\vec y}}&=\innerproduct{\sum_{j=1}^d\innerproduct{\widetilde x_j-\widetilde y_j}^2}^{\frac{1}{2}}, \\
		G_d\innerproduct{\widetilde{\vec x},\widetilde{\vec y}}&=\begin{cases} \dfrac{\imagunit}{4}H_0^{(1)}\innerproduct{k\rho\innerproduct{\widetilde{\vec x},\widetilde{\vec y}}}, & d=2, \\[2mm]
			\dfrac{e^{\imagunit k\rho\innerproduct{\widetilde{\vec x},\widetilde{\vec y}}}}{4\pi\rho\innerproduct{\widetilde{\vec x},\widetilde{\vec y}}}, & d=3. \end{cases}
	\end{align*}
	The square root is chosen so that $\Rea\rho>0$ and $\Ima\rho\geq0$. Here $H_0^{(1)}$ is the Hankel function of the first kind of order zero. Appendix \ref{app:green-truncation} gives the fundamental-solution identity and the kernel estimates needed below.
	\begin{theorem}
		\label{thm:truncated-optimal}
		Suppose that
		\begin{align}
			D_c\geq\bar l, \qquad k\geq k_0, \label{eq:green-absorption}
		\end{align}
		and that $k\sigma_0D_c$ is sufficiently large. Then \eqref{eq:truncated-weak} has a unique solution, and
		\begin{align}
			\enorm{u_L}_{B_L}\leq2^{\frac{d}{2}+2}C_\infty k\norm F_{V_L'}. \label{eq:truncated-stability}
		\end{align}
		Moreover, $\mu_L\eqsim k^{-1}$.
	\end{theorem}
	\begin{proof}
		Fix $w\in V_L$, and define $F_w\in V_L'$ by
		\begin{align*}
			\coproduct{F_w,\psi}=a_L\innerproduct{w,\psi}, \qquad \psi\in V_L.
		\end{align*}
		Let $\widetilde w$ be the reflection extension from Lemma \ref{lem:reflection-extension}. For $\psi\in H_0^1(B_{L+D_c})$, set
		\begin{align*}
			\coproduct{\widetilde F_w,\psi}=a_{B_{L+D_c}}\innerproduct{\widetilde w,\psi}.
		\end{align*}
		By \eqref{eq:reflection-pairing},
		\begin{align*}
			\snorm{\coproduct{\widetilde F_w,\psi}}\leq2^{\frac{d}{2}}\norm{F_w}_{V_L'}\enorm\psi_{B_{L+D_c}}.
		\end{align*}
		The Hahn--Banach theorem \cite[Theorem~5.6]{Folland} yields a norm-preserving extension $F_w^e$ of $\widetilde F_w$ to $V$. By Theorem \ref{thm:high-wellposed}, there exists $v\in V$ such that
		\begin{align}
			a\innerproduct{v,\psi}=\coproduct{F_w^e,\psi}, \qquad \psi\in V, \qquad \enorm v\leq2^{\frac{d}{2}}C_\infty k\norm{F_w}_{V_L'}. \label{eq:green-v-bound}
		\end{align}
		Set $z=v-\widetilde w$ on $B_{L+D_c}$. Since $F_w^e=\widetilde F_w$ on $H_0^1(B_{L+D_c})$,
		\begin{align}
			a_{B_{L+D_c}}\innerproduct{z,\psi}=0, \qquad \psi\in H_0^1(B_{L+D_c}). \label{eq:green-z-homogeneous}
		\end{align}
		Lemma \ref{lem:green-interior-decay} shows that, for $k\sigma_0D_c$ sufficiently large,
		\begin{align}
			\enorm z_{B_L}\leq2^{-\frac{d}{2}-1}\enorm z_{B_{L+D_c}}\leq\frac12\innerproduct{C_\infty k\norm{F_w}_{V_L'}+\enorm w_{B_L}}, \label{eq:green-z-bound}
		\end{align}
		where the second inequality uses \eqref{eq:reflection-extension} and \eqref{eq:green-v-bound}. Since $w=v-z$ on $B_L$, \eqref{eq:green-v-bound} and \eqref{eq:green-z-bound} give
		\begin{align*}
			\enorm w_{B_L}&\leq\enorm v+\enorm z_{B_L}\\
			&\leq\innerproduct{2^{\frac{d}{2}}+\frac12}C_\infty k\norm{F_w}_{V_L'}+\frac12\enorm w_{B_L}.
		\end{align*}
		Hence,
		\begin{align}
			\enorm w_{B_L}\leq2^{\frac{d}{2}+2}C_\infty k\norm{F_w}_{V_L'}. \label{eq:green-functional-bound}
		\end{align}
		Then \eqref{eq:green-functional-bound} gives $\mu_L\eqsim k^{-1}$.
	\end{proof}
	
	For $u,v\in V_L$, set
	\begin{align*}
		b_L\innerproduct{u,v}=\innerproduct{\mat A\grad u,\grad v}_{B_L}.
	\end{align*}
	\begin{theorem}\label{thm:cip-regularity}
		Let $g\in L^2(B_L)$, and let $w\in V_L$ solve $a_L\innerproduct{w,v}=\innerproduct{g,v}_{B_L}$ for all $v\in V_L$. Under the conditions of Theorem \ref{thm:truncated-optimal}, we have $w\in H^2(B_L)$ and
		\begin{align}
			\enorm w_{B_L}+k^{-1}\norm w_{H^2(B_L)}\leq C\norm g_{L^2(B_L)}. \label{eq:cip-regularity}
		\end{align}
	\end{theorem}
	\begin{proof}
		By Theorem \ref{thm:truncated-optimal} and $\norm v_{L^2(B_L)}\leq k^{-1}\enorm v_{B_L}$,
		\begin{align*}
			\enorm w_{B_L}\leq Ck\sup_{0\ne v\in V_L}\frac{\snorm{\innerproduct{g,v}_{B_L}}}{\enorm v_{B_L}}\leq C\norm g_{L^2(B_L)}.
		\end{align*}
		Since $\mat A,J$ are continuous and the right-angle corners of $B_L$ produce no singular exponents in the local Mellin-pencil analysis \cite{NicaiseSandig1994}, the proof of \cite[Theorem~6.3]{WangZheng} extends to $s=1$, giving $w\in H^2(B_L)$ with
		\begin{align*}
			k^{-1}\norm w_{H^2(B_L)}\leq C\norm g_{L^2(B_L)}.
		\end{align*}
		This proves \eqref{eq:cip-regularity}. For $b_L$, the Poincar\'e inequality and the same regularity argument without the mass term yield
		\begin{align*}
			\norm w_{H^1(B_L)}\leq C\norm g_{L^2(B_L)},\qquad \norm w_{H^2(B_L)}\leq C\norm g_{L^2(B_L)}.
		\end{align*}
	\end{proof}
	
	We next estimate the PML truncation error, following the boundary-lifting argument in \cite[Theorem~3.7]{LiWu}. We equip $H^{\frac{1}{2}}(\partial B_L)$ with the norm \cite[Chapter~3]{McLeanElliptic}
	\begin{align*}
		\norm g_{H^{\frac{1}{2}}(\partial B_L)}^2=\norm g_{L^2(\partial B_L)}^2+\int_{\partial B_L}\int_{\partial B_L}\frac{\snorm{g\innerproduct{\vec x}-g\innerproduct{\vec y}}^2}{\snorm{\vec x-\vec y}^{d}}\dd S_{\vec x}\dd S_{\vec y}.
	\end{align*}
	Here $H^{\frac{1}{2}}(\partial B_L)$ consists of all $g\in L^2(\partial B_L)$ for which this norm is finite.
	\begin{lemma}\label{lem:pml-lifting}
		Let $g\in H^{\frac12}(\partial B_L)$. Then there exists $\phi\in H^1(B_L)$ such that $\phi|_{\partial B_L}=g$ and
		\begin{align}
			\enorm\phi_{B_L}\leq C\innerproduct{\norm g_{H^{\frac12}(\partial B_L)}+k^{\frac12}\norm g_{L^2(\partial B_L)}}. \label{eq:pml-lifting}
		\end{align}
	\end{lemma}
	\begin{proof}
		By the trace theorem \cite[Theorem~3.37]{McLeanElliptic}, there exists $\psi\in H^1(B_L)$ such that $\psi|_{\partial B_L}=g$ and
		\begin{align*}
			\norm\psi_{H^1(B_L)}\leq C\norm g_{H^{\frac12}(\partial B_L)}.
		\end{align*}
		Since $k^{-1}<L$, let
		\begin{align*}
			S_k=\setprod{\vec x\in B_L:\operatorname{dist}(\vec x,\partial B_L)<k^{-1}}.
		\end{align*}
		We claim that
		\begin{align}
			\norm\psi_{L^2(S_k)}^2
			\leq\frac2k\norm g_{L^2(\partial B_L)}^2+\frac1{k^2}\norm{\grad\psi}_{L^2(B_L)}^2. \label{eq:pml-strip-bound}
		\end{align}
		For $j=1,\ldots,d$, let $F_j^\pm=\setprod{\vec y\in\partial B_L:y_j=\pm L}$, and let $\vec e_j$ denote the $j$th coordinate unit vector. $S_k$ is covered by the images of $\vec x=\vec y-t\vec n$, where $\vec y\in F_j^\pm$, $\vec n=\pm\vec e_j$ is the outward normal, and $0<t<k^{-1}$. For almost every $\vec y$,
		\begin{align*}
			\psi(\vec y-t\vec n)&=g(\vec y)-\int_0^t\partial_{\vec n}\psi(\vec y-s\vec n)\dd s,\\
			\snorm{\psi(\vec y-t\vec n)}^2&\leq2\snorm{g(\vec y)}^2+2t\int_0^t\snorm{\partial_{\vec n}\psi(\vec y-s\vec n)}^2\dd s.
		\end{align*}
		Integrating in $t$ gives
		\begin{align*}
			\int_0^{k^{-1}}\snorm{\psi(\vec y-t\vec n)}^2\dd t
			\leq\frac2k\snorm{g(\vec y)}^2+\frac1{k^2}\int_0^{k^{-1}}\snorm{\partial_{\vec n}\psi(\vec y-s\vec n)}^2\dd s.
		\end{align*}
		Integrating the pointwise bound over each face and summing over the $2d$ faces $x_j=\pm L$,
		\begin{align*}
			\norm\psi_{L^2(S_k)}^2
			&\leq\sum_{j=1}^d\sum_{\pm}\int_{F_j^\pm}\int_0^{k^{-1}}\snorm{\psi\innerproduct{\vec y-t\vec n}}^2\dd t\dd S_{\vec y}\\
			&\leq\frac2k\sum_{j=1}^d\sum_{\pm}\int_{F_j^\pm}\snorm{g\innerproduct{\vec y}}^2\dd S_{\vec y}\\
			&\qquad+\frac1{k^2}\sum_{j=1}^d\sum_{\pm}\int_{F_j^\pm}\int_0^{k^{-1}}\snorm{\partial_{\vec n}\psi\innerproduct{\vec y-s\vec n}}^2\dd s\dd S_{\vec y}\\
			&\leq\frac2k\norm g_{L^2(\partial B_L)}^2+\frac1{k^2}\sum_{j=1}^d\norm{\partial_{x_j}\psi}_{L^2(B_L)}^2\\
			&=\frac2k\norm g_{L^2(\partial B_L)}^2+\frac1{k^2}\norm{\grad\psi}_{L^2(B_L)}^2,
		\end{align*}
		which proves \eqref{eq:pml-strip-bound}.
		
		Let
		\begin{align*}
			\chi_k(\vec x)=\max\setprod{1-k\operatorname{dist}(\vec x,\partial B_L),0},\qquad \phi=\chi_k\psi.
		\end{align*}
		Since the distance function is $1$-Lipschitz, $\chi_k\in W^{1,\infty}(B_L)$, $0\leq\chi_k\leq1$ and $\snorm{\grad\chi_k}\leq k$; moreover, $\chi_k=1$ on $\partial B_L$ and $\chi_k=0$ on $B_L\setminus S_k$, so $\phi\in H^1(B_L)$ and $\phi|_{\partial B_L}=g$. By the product rule,
		\begin{align*}
			\enorm\phi_{B_L}^2
			&=\int_{S_k}\innerproduct{\snorm{\chi_k\grad\psi+\psi\grad\chi_k}^2+k^2\snorm{\chi_k\psi}^2}\dd\vec x\\
			&\leq2\int_{S_k}\chi_k^2\snorm{\grad\psi}^2\dd\vec x+2\int_{S_k}\snorm\psi^2\snorm{\grad\chi_k}^2\dd\vec x+k^2\int_{S_k}\chi_k^2\snorm\psi^2\dd\vec x\\
			&\leq2\norm{\grad\psi}_{L^2(B_L)}^2+3k^2\norm\psi_{L^2(S_k)}^2\\
			&\leq5\norm{\grad\psi}_{L^2(B_L)}^2+6k\norm g_{L^2(\partial B_L)}^2\\
			&\leq C\innerproduct{\norm g_{H^{\frac12}(\partial B_L)}^2+k\norm g_{L^2(\partial B_L)}^2}.
		\end{align*}
		This proves \eqref{eq:pml-lifting}.
	\end{proof}
	
	\begin{theorem}
		\label{thm:pml-truncation-error}
		For $f\in L^2(\R^d)$ compactly supported in $B_l$, let $u$ denote the solution of \eqref{eq:helm-original}--\eqref{eq:sommerfeld}. Under the conditions of Theorem \ref{thm:truncated-optimal}, let $\widetilde u\in V$ solve \eqref{eq:UPMLstrong} and $u_L\in V_L$ solve \eqref{eq:truncated-weak} with $\coproduct{F,v}=\innerproduct{f,v}_{B_l}$. Set
		\begin{align*}
			\beta_2=\frac{1}{6\sqrt2},\qquad \beta_3=\frac{1}{3\sqrt3}.
		\end{align*}
		Then
		\begin{align}
			\enorm{\widetilde u-u_L}_{B_L}\leq Ck^{\frac d2}e^{-\beta_dk\sigma_0D_c}\norm f_{L^2(B_l)}. \label{eq:pml-truncation-error}
		\end{align}
		The constant $C$ is independent of $k$ and $f$.
	\end{theorem}
	\begin{proof}
		The Green representation reads
		\begin{align}
			\widetilde u(\vec x)=\int_{B_l}G_d(\widetilde{\vec x},\vec y)f(\vec y)\dd\vec y. \label{eq:pml-source-representation}
		\end{align}
		Set $w=\widetilde u-u_L$ in $B_L$ and $g=\widetilde u|_{\partial B_L}$. Since $f$ is supported in $B_l$, the two equations read
		\begin{align*}
			a_L\innerproduct{\widetilde u,v}&=\innerproduct{f,v}_{B_l},\qquad v\in V_L;\\
			a_L\innerproduct{u_L,v}&=\innerproduct{f,v}_{B_l},\qquad v\in V_L.
		\end{align*}
		Subtracting them and using $u_L|_{\partial B_L}=0$ yields
		\begin{align*}
			a_L\innerproduct{w,v}=0,\qquad \forall v\in V_L;\qquad w|_{\partial B_L}=g.
		\end{align*}
		By Lemma \ref{lem:pml-lifting}, there exists $\phi\in H^1(B_L)$ with $\phi|_{\partial B_L}=g$ and
		\begin{align*}
			\enorm\phi_{B_L}\leq C\innerproduct{\norm g_{H^{\frac12}(\partial B_L)}+k^{\frac12}\norm g_{L^2(\partial B_L)}}.
		\end{align*}
		Since $w-\phi\in V_L$ and $a_L\innerproduct{w-\phi,v}=-a_L\innerproduct{\phi,v}$ for all $v\in V_L$, Theorem \ref{thm:truncated-optimal} and the coefficient bounds \eqref{eq:truncated-coefficient-bounds} give
		\begin{align*}
			\enorm{w-\phi}_{B_L}
			&\leq Ck\sup_{0\ne v\in V_L}\frac{\snorm{a_L\innerproduct{\phi,v}}}{\enorm v_{B_L}}
			\leq Ck\enorm\phi_{B_L}.
		\end{align*}
		Hence
		\begin{align*}
			\enorm w_{B_L}\leq\enorm\phi_{B_L}+\enorm{w-\phi}_{B_L}
			\leq Ck\innerproduct{\norm g_{H^{\frac12}(\partial B_L)}+k^{\frac12}\norm g_{L^2(\partial B_L)}}.
		\end{align*}
		
		It remains to estimate $g$. For $\vec x\in\partial B_L$ and $\vec y\in B_l$, choose $j$ with $\snorm{x_j}=L$. Since $\sigma=\sigma_0$ for $\snorm t\geq\bar l$ and $D_c\geq\bar l>l$,
		\begin{align*}
			\snorm{\vec x-\vec y}\geq L-l\geq D_c,\quad \snorm{\vec x-\vec y}\leq\sqrt d\innerproduct{L+l}\leq3\sqrt dD_c,\quad \snorm{\int_{y_j}^{x_j}\sigma\innerproduct t\dd t}\geq\sigma_0D_c.
		\end{align*}
		Consequently, \eqref{eq:rho-imaginary} gives
		\begin{align*}
			\Ima\rho\innerproduct{\widetilde{\vec x},\vec y}\geq\frac{\sigma_0D_c^2}{\snorm{\vec x-\vec y}}\geq\frac{\sigma_0D_c}{3\sqrt d}.
		\end{align*}
		As in the proof of Lemma \ref{lem:green-interior-decay}, take $k_0\bar l\geq2$, so $kD_c\geq2$. Thus, for $m=0,1$, estimates \eqref{eq:G2-pointwise}--\eqref{eq:G3-pointwise} give
		\begin{align*}
			\snorm{\grad_{\vec x}^{m}G_d\innerproduct{\widetilde{\vec x},\vec y}}\leq Ck^{m+\frac{d-3}{2}}e^{-\beta_dk\sigma_0D_c},\qquad \vec x\in\partial B_L,\quad\vec y\in B_l.
		\end{align*}
		Differentiating \eqref{eq:pml-source-representation} on $\partial B_L$ gives
		\begin{align*}
			\norm g_{L^2(\partial B_L)}&\leq Ck^{\frac{d-3}{2}}e^{-\beta_dk\sigma_0D_c}\norm f_{L^2(B_l)},\\
			\norm g_{H^1(\partial B_L)}&\leq Ck^{\frac{d-1}{2}}e^{-\beta_dk\sigma_0D_c}\norm f_{L^2(B_l)}.
		\end{align*}
		Interpolation on the Lipschitz boundary \cite[Chapter~3]{McLeanElliptic} gives
		\begin{align}
			\norm g_{H^{\frac12}(\partial B_L)}&\leq C\norm g_{L^2(\partial B_L)}^{\frac12}\norm g_{H^1(\partial B_L)}^{\frac12}\notag\\
			&\leq Ck^{\frac{d-2}{2}}e^{-\beta_dk\sigma_0D_c}\norm f_{L^2(B_l)}. \label{eq:pml-boundary-trace}
		\end{align}
		Substituting the two trace bounds into the estimate for $w$ yields
		\begin{align*}
			\enorm w_{B_L}&\leq Ck\innerproduct{k^{\frac{d-2}{2}}+k^{\frac12}k^{\frac{d-3}{2}}}e^{-\beta_dk\sigma_0D_c}\norm f_{L^2(B_l)}
			\leq Ck^{\frac d2}e^{-\beta_dk\sigma_0D_c}\norm f_{L^2(B_l)},
		\end{align*}
		which proves \eqref{eq:pml-truncation-error}.
	\end{proof}
	
	\section{Preasymptotic error analysis of the CIP-FEM}
	\label{sec:cip-fem}
	Throughout this section, we work under the conditions of Theorem \ref{thm:truncated-optimal}, namely the regularity conditions on $\sigma$ from Section \ref{sec:model} and $D_c\geq\bar l$, $k\geq k_0$ with $k\sigma_0D_c$ sufficiently large. Consequently, the stability estimate $\mu_L\eqsim k^{-1}$ holds. We analyze the linear CIP-FEM following \cite[Section~4]{LiWu}. Constants are independent of $k$, $h$, and $f$, but may depend on the fixed PML parameters, mesh regularity, and the penalty bound. Here $f\in L^2(\R^d)$ is compactly supported in $B_l$.
	
	Let $\mathcal T_h$ be a shape-regular, quasi-uniform mesh of $B_L$. Let $V_h\subset V_L$ be the corresponding continuous $P_1$ or $Q_1$ finite element space. 
	For each element $K\in\mathcal T_h$, let $h_K=\operatorname{diam}K$, and set $h=\max_{K\in\mathcal T_h}h_K$. Denote by $\mathcal E_h^I$ the set of interior edges or faces. For each $e\in\mathcal E_h^I$, let $h_e=\operatorname{diam}e$ and fix a unit normal $\vec n_e$.
	
	Let the penalty parameters satisfy $|\gamma_e|\leq C_\gamma$ and $\Rea\innerproduct{e^{-\imagunit\vartheta_d}\gamma_e}\geq-\gamma_0$, where $C_\gamma,\gamma_0>0$ are independent of $h$ and $k$. For $d=2$, this reduces to the condition in \cite[(6.11)]{WangZheng}. Set
	\begin{align*}
		b_h=b_L+\mathcal J_h,\qquad a_h\innerproduct{w,v}=b_h\innerproduct{w,v}-k^2\innerproduct{Jw,v}_{B_L},
	\end{align*}
	where
	\begin{align*}
		\mathcal J_h\innerproduct{w,v}=\sum_{e\in\mathcal E_h^I}\gamma_eh_e \innerproduct{\jump{\grad w\cdot\vec n_e},\jump{\grad v\cdot\vec n_e}}_e.
	\end{align*}
	Here $\llbracket\cdot\rrbracket$ denotes the jump of the normal derivative across $e$ \cite[Section~2]{Wu2014}. The CIP-FEM seeks $u_{L,h}\in V_h$ such that
	\begin{align}
		a_h\innerproduct{u_{L,h},v_h}=\innerproduct{f,v_h}_{B_l},\qquad v_h\in V_h. \label{eq:cip-discrete}
	\end{align}
	The case $\gamma_e=0$ gives the standard FEM.
	
	For $w\in V_L$, define the elliptic projections $P_h^\pm w\in V_h$ by
	\begin{align*}
		b_h\innerproduct{P_h^+w,v_h}=b_h\innerproduct{w,v_h},\qquad b_h\innerproduct{v_h,P_h^-w}=b_h\innerproduct{v_h,w},\qquad v_h\in V_h.
	\end{align*}
	\begin{lemma}\label{lem:cip-projection}
		There exists a constant \( \gamma_0 >0 \)  such that if $\Rea\innerproduct{e^{-\imagunit\vartheta_d}\gamma_e}\geq-\gamma_0$, we have
		\begin{align}
			\norm{w-P_h^\pm w}_{L^2(B_L)}+h\snorm{w-P_h^\pm w}_{H^1(B_L)} +h^2\snorm{w-P_h^\pm w}_{H^2(\mathcal T_h)} \lesssim h^2\snorm w_{H^2(B_L)}, \label{eq:cip-projection}
		\end{align}
		where $\snorm v_{H^2(\mathcal T_h)}^2=\sum_K\snorm v_{H^2(K)}^2$.
	\end{lemma}
	\begin{proof}
		Following \cite[Lemma~4.2]{LiWu}, the elementwise trace inequality and \eqref{eq:Aelliptic} give
		\begin{align*}
			\snorm{b_h\innerproduct{w,v}}&\lesssim\innerproduct{\snorm w_{H^1(B_L)}+h\snorm w_{H^2(\mathcal T_h)}}\innerproduct{\snorm v_{H^1(B_L)}+h\snorm v_{H^2(\mathcal T_h)}},\\
			\Rea\innerproduct{e^{-\imagunit\vartheta_d}b_h\innerproduct{v_h,v_h}} &\geq a_*\snorm{v_h}_{H^1(B_L)}^2+\sum_{e\in\mathcal E_h^I}\Rea\innerproduct{e^{-\imagunit\vartheta_d}\gamma_e}h_e\norm{\jump{\grad v_h\cdot\vec n_e}}_{L^2(e)}^2\\
			&\geq a_*\snorm{v_h}_{H^1(B_L)}^2-\gamma_0\sum_{e\in\mathcal E_h^I}h_e\norm{\jump{\grad v_h\cdot\vec n_e}}_{L^2(e)}^2\\
			&\geq\innerproduct{a_*-\gamma_0C}\snorm{v_h}_{H^1(B_L)}^2\geq\frac{a_*}2\snorm{v_h}_{H^1(B_L)}^2,
		\end{align*}
		where the last step uses the elementwise trace inequality with a constant $C$ independent of $h$ and $k$, and $\gamma_0$ is chosen so that $\gamma_0C\leq a_*/2$.
		
		Let $I_h$ denote the Lagrange interpolation operator. The standard interpolation estimate \cite[Chapter~3]{Ciarlet} is
		\begin{align*}
			\norm{w-I_hw}_{L^2(B_L)}+h\snorm{w-I_hw}_{H^1(B_L)}+h^2\snorm{w-I_hw}_{H^2(\mathcal T_h)}\lesssim h^2\snorm w_{H^2(B_L)}.
		\end{align*}
		For $P_h^-$, set $\xi_h=I_hw-P_h^-w$. The inverse estimate gives
		\begin{align*}
			h\snorm{\xi_h}_{H^2(\mathcal T_h)}\lesssim\snorm{\xi_h}_{H^1(B_L)}.
		\end{align*}
		Coercivity, projection orthogonality, and the preceding bounds give
		\begin{align*}
			\frac{a_*}{2}\snorm{\xi_h}_{H^1(B_L)}^2&\leq\snorm{b_h\innerproduct{\xi_h,\xi_h}} =\snorm{b_h\innerproduct{\xi_h,I_hw-w}}\lesssim h\snorm{\xi_h}_{H^1(B_L)}\snorm w_{H^2(B_L)}.
		\end{align*}
		Hence, $\snorm{w-P_h^-w}_{H^1(B_L)}\lesssim h\snorm w_{H^2(B_L)}$. For $e=w-P_h^-w$, the interpolation and inverse estimates give
		\begin{align*}
			\snorm e_{H^2(\mathcal T_h)}&\leq\snorm{w-I_hw}_{H^2(\mathcal T_h)}+\snorm{\xi_h}_{H^2(\mathcal T_h)}\\
			&\lesssim\snorm w_{H^2(B_L)}+h^{-1}\snorm{\xi_h}_{H^1(B_L)}\lesssim\snorm w_{H^2(B_L)}.
		\end{align*}
		Let $z\in V_L$ solve $b_L\innerproduct{z,v}=\innerproduct{e,v}_{B_L}$. By Theorem \ref{thm:cip-regularity}, $\norm z_{H^2(B_L)}\lesssim\norm e_{L^2(B_L)}$. Orthogonality gives
		\begin{align*}
			\norm e_{L^2(B_L)}^2&=b_h\innerproduct{z,e}=b_h\innerproduct{z-I_hz,e}\\& \lesssim h^2\norm z_{H^2(B_L)}\snorm w_{H^2(B_L)}\lesssim h^2\norm e_{L^2(B_L)}\snorm w_{H^2(B_L)}.
		\end{align*}
		The estimate for $P_h^+$ follows by the same argument with the adjoint problem for $b_L$.
	\end{proof}
	
	\begin{theorem}\label{thm:cip-error}
		\( u_{L,h} \) and \( u_L \) are solutions to \eqref{eq:cip-discrete} and \eqref{eq:truncated-weak} with \( F=f \), respectively. There exist $c_0,C>0$, independent of $k,h,f$, such that, if $k^3h^2\leq c_0$, we have
		\begin{align}
			&\enorm{u_L-u_{L,h}}_{B_L}\leq C\innerproduct{kh+k^3h^2}\norm f_{L^2(B_l)}, \label{eq:cip-error-energy}\\ &\norm{u_L-u_{L,h}}_{L^2(B_L)}\leq Ck^2h^2\norm f_{L^2(B_l)}. \label{eq:cip-error-l2}
		\end{align}
		Moreover, $\enorm{u_{L,h}}_{B_L}\leq C\norm f_{L^2(B_l)}$.
	\end{theorem}
	\begin{proof}
		We follow the modified duality argument of \cite[Theorem~4.4]{LiWu}. Let $u_{L,h}$ be the discrete solution and set $e_h=u_L-u_{L,h}$. Since $u_L\in H^2(B_L)$, we have $\mathcal J_h\innerproduct{u_L,v}=0$ and hence $a_h\innerproduct{e_h,v_h}=0$. Let $z\in V_L$ solve $a_L\innerproduct{v,z}=\innerproduct{v,e_h}_{B_L}$. Theorem \ref{thm:cip-regularity} gives
		\begin{align*}
			\norm{u_L}_{H^2(B_L)}\lesssim k\norm f_{L^2(B_l)},\qquad \norm z_{H^2(B_L)}\lesssim k\norm{e_h}_{L^2(B_L)}.
		\end{align*}
		Using projection orthogonality, $\mathcal J_h\innerproduct{v,z}=0$, and Lemma \ref{lem:cip-projection}, we obtain
		\begin{align*}
			\norm{e_h}_{L^2(B_L)}^2&=a_h\innerproduct{e_h,z-P_h^-z}\\
			&=b_h\innerproduct{u_L-I_hu_L,z-P_h^-z}-k^2\innerproduct{Je_h,z-P_h^-z}_{B_L}\\
			&\lesssim h^2\norm{u_L}_{H^2(B_L)}\norm z_{H^2(B_L)}+k^2h^2\norm{e_h}_{L^2(B_L)}\norm z_{H^2(B_L)}\\
			&\leq Ck^2h^2\norm f_{L^2(B_l)}\norm{e_h}_{L^2(B_L)}+Ck^3h^2\norm{e_h}_{L^2(B_L)}^2.
		\end{align*}
		Choose $c_0$ so that $Cc_0\leq\frac{1}{2}$ and absorb the last term to obtain the $L^2(B_L)$ estimate in \eqref{eq:cip-error-l2}. For $\eta_h=u_{L,h}-P_h^+u_L$, Galerkin orthogonality and the projection property give
		\begin{align*}
			\frac{a_*}{2}\snorm{\eta_h}_{H^1(B_L)}^2&\leq\snorm{b_h\innerproduct{\eta_h,\eta_h}} =k^2\snorm{\innerproduct{Je_h,\eta_h}_{B_L}}\\
			&\lesssim k^2\norm{e_h}_{L^2(B_L)}\innerproduct{\norm{e_h}_{L^2(B_L)}+\norm{u_L-P_h^+u_L}_{L^2(B_L)}} \\
			&\lesssim k^6h^4\norm f_{L^2(B_l)}^2.
		\end{align*}
		This yields $\snorm{\eta_h}_{H^1(B_L)}\lesssim k^3h^2\norm f_{L^2(B_l)}$. Using this bound together with \eqref{eq:cip-error-l2} and Lemma \ref{lem:cip-projection}, we obtain
		\begin{align*}
			\enorm{e_h}_{B_L}&\leq\snorm{u_L-P_h^+u_L}_{H^1(B_L)}+\snorm{\eta_h}_{H^1(B_L)}+k\norm{e_h}_{L^2(B_L)}\\
			&\lesssim\innerproduct{kh+k^3h^2}\norm f_{L^2(B_l)}.
		\end{align*}
		Finally, $k^3h^2\leq c_0$ and $k\geq k_0>1$ imply $kh\lesssim1$, so $\enorm{u_{L,h}}_{B_L}\leq\enorm{u_L}_{B_L}+\enorm{e_h}_{B_L}\lesssim\norm f_{L^2(B_l)}$.
	\end{proof}
	
	\begin{theorem}
		\label{thm:helmholtz-cip-error}
		Let $\widetilde u$ and $u_{L,h}$ be solutions to \eqref{eq:UPMLstrong} and \eqref{eq:cip-discrete}, respectively. Under the conditions of Theorems \ref{thm:pml-truncation-error} and \ref{thm:cip-error},  we have
		\begin{align}
			\enorm{\widetilde u-u_{L,h}}_{B_L}\leq C\innerproduct{kh+k^3h^2+k^{\frac d2}e^{-\beta_dk\sigma_0D_c}}\norm f_{L^2(B_l)}. \label{eq:helmholtz-cip-error}
		\end{align}
		Here $\beta_d$ is as in Theorem \ref{thm:pml-truncation-error}, and $C$ is independent of $k$, $h$, and $f$. The estimate also holds for the standard FEM ($\gamma_e=0$).
	\end{theorem}
	\begin{proof}
		We combine the truncation and discretization estimates as in \cite[Remark~4.5(v)]{LiWu}. By the triangle inequality,
		\begin{align*}
			\enorm{\widetilde u-u_{L,h}}_{B_L}
			\leq\enorm{\widetilde u-u_L}_{B_L}+\enorm{u_L-u_{L,h}}_{B_L}.
		\end{align*}
		Setting $\gamma_e=0$ yields the standard FEM.
	\end{proof}
	
	\section{Numerical tests}\label{sec:numerics}
	In this section, we solve the two-dimensional Helmholtz problem \eqref{eq:helm-original}--\eqref{eq:sommerfeld}. The problem is truncated by the UPML, and the resulting problem is discretized by the standard FEM ($\gamma_e=0$) and by the CIP-FEM \eqref{eq:cip-discrete} on the same mesh. It should be noted that, in fact, a large number of numerical simulations of UPML have been conducted over the past three decades. Therefore, here we merely present a simple demonstration of the effectiveness of CIP-FEM in eliminating the pollution errors of UPML problems. All computations are performed with the finite element library MFEM \cite{mfem-web}.
	
	We set $B_l=(-0.5,0.5)^2$ and $B_L=(-1,1)^2$, so that $l=0.5$, $\bar l=0.8$, $D=0.3$, and $D_c=L-\bar l=0.2$. 	The absorption function is \( \sigma(t) = \sigma_0\innerproduct{\dfrac{\snorm t-l}{\bar l-l}}^2 (l<\snorm t<\bar l) \) with \( \sigma_0=40 \).	For the CIP-FEM we use the constant penalty parameter $\gamma_e=-1/24$.
	
	Let $r=\snorm{\vec x}$ and $r_0=0.25$. We choose the source term $f$ so that the exact solution of \eqref{eq:helm-original}--\eqref{eq:sommerfeld} is
	\begin{align}
		u(\vec x)=\frac{\imagunit}{4}H_0^{(1)}(kr)\,
		\begin{cases}
			\dfrac{1}{1+\exp\innerproduct{\dfrac{1-2r/r_0}{(r/r_0)\innerproduct{1-r/r_0}}}}, & r<r_0,\\[3mm]
			1, & r\geq r_0,
		\end{cases}\label{eq:num-exact}
	\end{align}
	where $H_0^{(1)}$ denotes the Hankel function of the first kind of order zero.

	We measure the error in the relative energy norm
	\begin{align*}
		E_h=\frac{\enorm{u-u_h}_{B_l}}{\enorm u_{B_l}},
	\end{align*}
	where $u_h\in V_h$ is the FEM or the CIP-FEM solution and $u$ is the exact solution \eqref{eq:num-exact}. For comparison we also report the error $E_h^{\rm int}$ of the nodal interpolant $I_hu\in V_h$.
	
\begin{figure}
		\centering
		\includegraphics[width=0.45\linewidth]{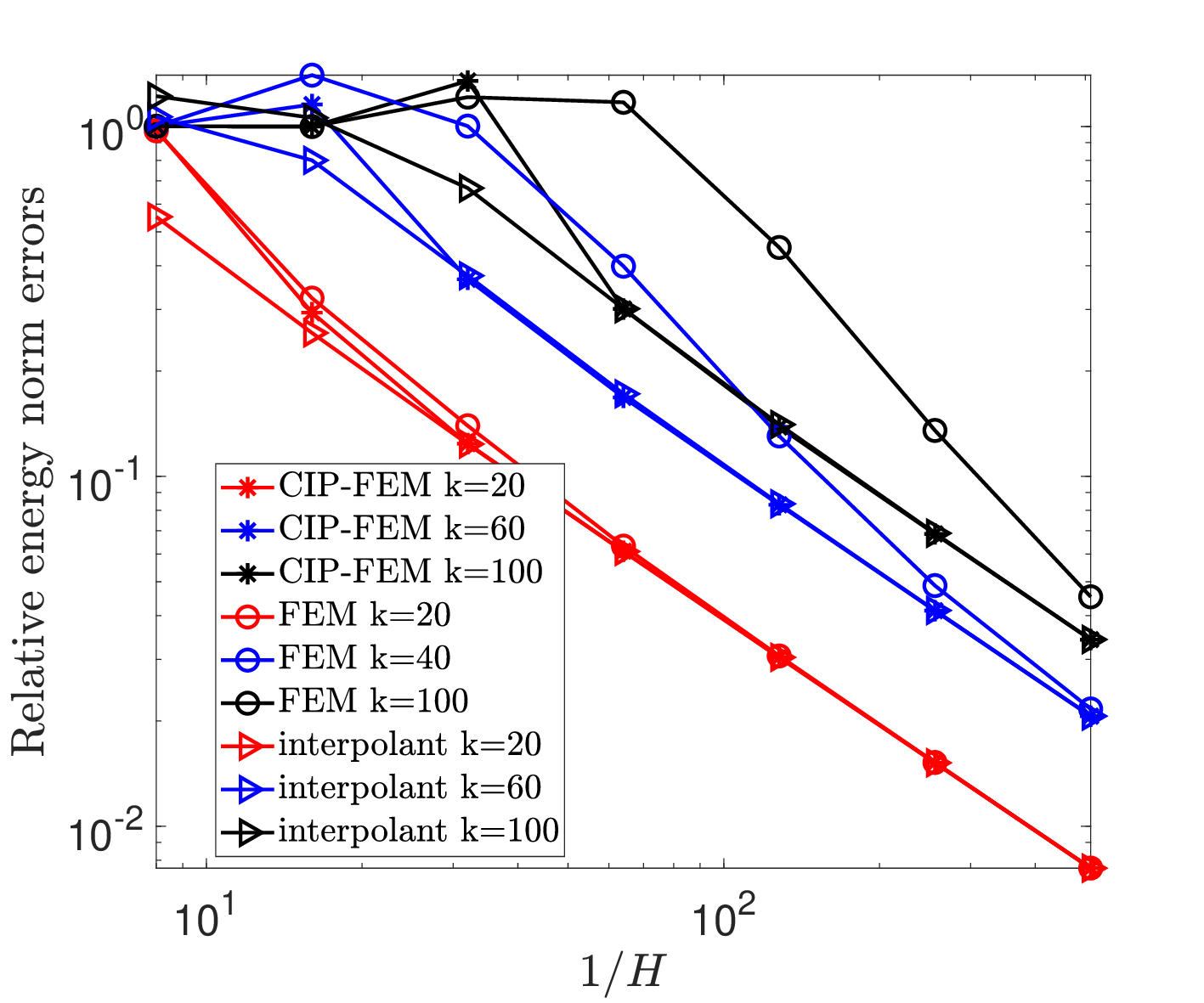}
		\hspace{0.02\linewidth}
		\includegraphics[width=0.45\linewidth]{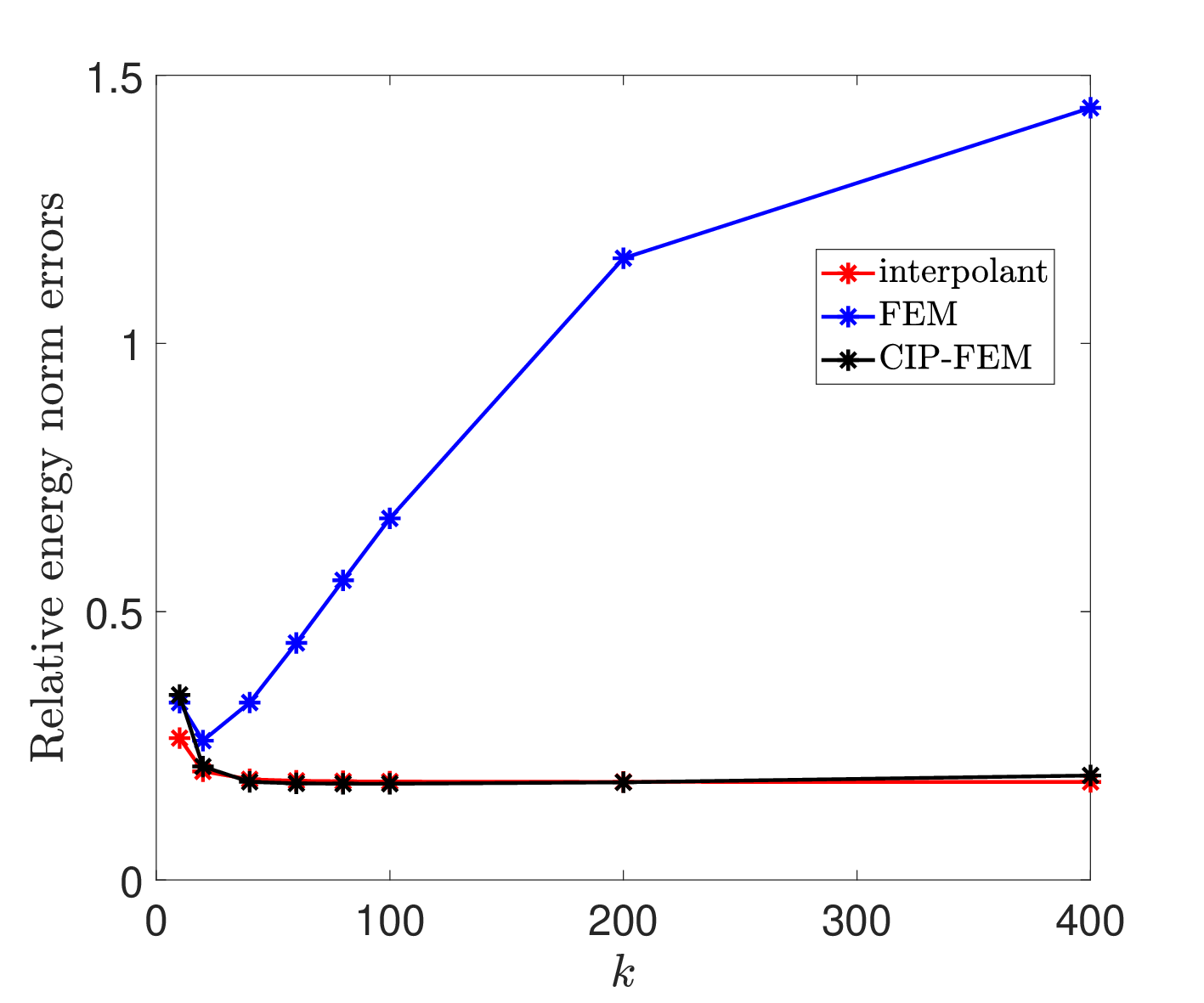}
		\caption{Relative energy-norm errors $E_h$ of the FEM, the CIP-FEM ($\gamma_e=-1/24$), and the nodal interpolant $I_hu$ for $l=0.5$, $\bar l=0.8$, $L=1$, and $\sigma_0=40$. Left: $E_h$ versus \( 1/h \) for $k=20,60,100$. Right: $E_h$ versus $k$ for fixed $kh=1$.}
		\label{fig:errors}
	\end{figure}

	In the left of Figure~\ref{fig:errors} we plot $E_h$ for $k=20,60,100$. For $k=20$ the CIP-FEM curve is almost indistinguishable from that of the FEM, whereas for $k=100$ it stays close to the interpolation curve down to meshes on which the FEM error is already dominated by pollution. Both methods are covered by the same bound \eqref{eq:cip-error-energy}; the numerical results indicate that the constant multiplying $k^3h^2$ is substantially smaller for the CIP-FEM, so that its error follows the best-approximation curve over a much wider range of mesh sizes, in particular at large $k$.
	
	In the right of Figure~\ref{fig:errors} we fix $kh=1$, i.e., $h=1/k$, and let $k$ increase from $10$ to $400$. In this scaling $kh$ is constant and $k^3h^2=k$, so that the standard FEM is expected to show an error that grows with $k$, whereas the CIP-FEM should remain at the level of the interpolation error. The FEM error departs from the interpolation error once $k$ exceeds about $20$ and then grows steadily, while the CIP-FEM error stays close to the interpolation error up to $k=400$.
	
	\appendix
	
	\section{Properties of the coefficients}
	\label{app:coefficients}
	
	\begin{lemma}
		\label{lem:coefficients}
		The coefficients satisfy $b_j,e_j\in W^{1,\infty}(\R^d)$ and $w_j\in L^\infty(\R^d)\cap C^{0,\frac{1}{2}}(\R^d)$. In addition,
		\begin{align}
			\frac2{M^4}\sigma_j\leq d_j,\quad \snorm{e_j}\leq Cd_j,\quad c_j+\sigma_0d_j\geq M^{-2},\quad j=1,\ldots,d. \label{eq:coeff-bounds}
		\end{align}
	\end{lemma}
	
	\begin{proof}
		Since $b_j\in W^{1,\infty}$, we have $e_j=b_j^2-1\in W^{1,\infty}$. The bounds $1\leq\rho_j\leq1+\zeta\sigma_0$ and $1\leq R_j\leq M^2$ give
		\begin{align*}
			\frac2{M^4}\sigma_j\leq\frac{2\rho_j\sigma_j}{R_j^2}\leq2\innerproduct{1+\zeta\sigma_0}\sigma_j,
		\end{align*}
		which proves the first inequality. Moreover,
		\begin{align*}
			e_j &=\frac{1-\innerproduct{1+\innerproduct{\zeta+\imagunit}\sigma_j}^2}{\innerproduct{1+\innerproduct{\zeta+\imagunit}\sigma_j}^2} =-\frac{\innerproduct{\zeta+\imagunit}\sigma_j\innerproduct{2+\innerproduct{\zeta+\imagunit}\sigma_j}}{\alpha_j^2}.
		\end{align*}
		Thus, we obtain $\snorm{e_j}\leq C\sigma_j$. The second inequality follows from the first. Finally,
		\begin{align*}
			\innerproduct{c_j+\sigma_0d_j}-R_j^{-1} &=\frac{\rho_j^2-\sigma_j^2+2\sigma_0\rho_j\sigma_j-\innerproduct{\rho_j^2+\sigma_j^2}}{R_j^2} =\frac{2\sigma_j\innerproduct{\sigma_0\rho_j-\sigma_j}}{R_j^2}\geq0.
		\end{align*}
		This proves the third inequality.
	\end{proof}
	
	\begin{lemma}\label{lm:c0}
		The functions $w_j$ defined in Section~\ref{sec:model} satisfy
		\begin{align}
			\norm{w_j}_{L^\infty}\leq \sqrt{2\innerproduct{1+\zeta\sigma_0}\sigma_0},\qquad \midproduct{w_j}_{C^{0,\frac{1}{2}}}\leq C_0, \label{eq:wjholder}
		\end{align}
		where $C_0$ is independent of $\hbar$ and $\delta$.
	\end{lemma}
	\begin{proof}
		The function $w_j$ depends only on $x_j$. Since $0\leq\sigma\leq\sigma_0$, we have
		\begin{align*}
			0\leq w_j\innerproduct{\vec x} =\frac{\sqrt{2\innerproduct{1+\zeta\sigma\innerproduct{x_j}}\sigma\innerproduct{x_j}}} {\innerproduct{1+\zeta\sigma\innerproduct{x_j}}^2+\sigma\innerproduct{x_j}^2} \leq\sqrt{2\innerproduct{1+\zeta\sigma_0}\sigma_0}.
		\end{align*}
		For any $\vec x,\vec y\in\R^d$, it follows that
		\begin{align*}
			\snorm{w_j\innerproduct{\vec x}-w_j\innerproduct{\vec y}} &\leq\snorm{w_j\innerproduct{\vec x}^2-w_j\innerproduct{\vec y}^2}^{\frac{1}{2}}\\
			&=\snorm{ \frac{2\innerproduct{1+\zeta\sigma\innerproduct{x_j}}\sigma\innerproduct{x_j}} {\innerproduct{\innerproduct{1+\zeta\sigma\innerproduct{x_j}}^2+\sigma\innerproduct{x_j}^2}^2} -\frac{2\innerproduct{1+\zeta\sigma\innerproduct{y_j}}\sigma\innerproduct{y_j}} {\innerproduct{\innerproduct{1+\zeta\sigma\innerproduct{y_j}}^2+\sigma\innerproduct{y_j}^2}^2} }^{\frac{1}{2}}\\
			&=\snorm{\int_{\sigma\innerproduct{y_j}}^{\sigma\innerproduct{x_j}} \frac{2-2\innerproduct{1+\zeta^2}s^2\innerproduct{3+2\zeta s}} {\innerproduct{\innerproduct{1+\zeta s}^2+s^2}^3}\dd s}^{\frac{1}{2}}\\
			&\leq C\snorm{\sigma\innerproduct{x_j}-\sigma\innerproduct{y_j}}^{\frac{1}{2}}.
		\end{align*}
		Consequently,
		\begin{align*}
			\snorm{\sigma\innerproduct{x_j}-\sigma\innerproduct{y_j}}^{\frac{1}{2}}\leq \norm{\sigma'}_{L^\infty}^{\frac{1}{2}}\snorm{x_j-y_j}^{\frac{1}{2}}\leq C\Theta\innerproduct{\sigma_0}^{\frac{1}{2}}\snorm{\vec x-\vec y}^{\frac{1}{2}}.
		\end{align*}
		From the definition of $\midproduct{w_j}_{C^{0,\frac{1}{2}}}$, we obtain
		\begin{align*}
			\midproduct{w_j}_{C^{0,\frac{1}{2}}} =\sup_{\vec x\ne\vec y} \frac{\snorm{w_j\innerproduct{\vec x}-w_j\innerproduct{\vec y}}} {\snorm{\vec x-\vec y}^{\frac{1}{2}}} \leq C\Theta\innerproduct{\sigma_0}^{\frac{1}{2}}\le C_0.
		\end{align*}
		This proves \eqref{eq:wjholder}.
	\end{proof}
	
	\section{\texorpdfstring{Fourier multipliers and the smoothing operator}{Fourier multipliers and the smoothing operator}}
	\label{app:fourier}
	\begin{lemma}
		\label{lem:fourier-operators}
		Fix $0<\hbar\le1$. The operators $\D_j$ and $P_{0,\hbar}$ are self-adjoint on the domains
		\begin{align*}
			\Dom\innerproduct{\D_j} &=\setprod{u\in L^2:\xi_j\widehat u\in L^2} =\setprod{u\in L^2:\partial_ju\in L^2},\\
			\Dom\innerproduct{P_{0,\hbar}} &=\setprod{u\in L^2:p_0\widehat u\in L^2}=H^2,
		\end{align*}
		respectively, where $\partial_j u$ is understood in the sense of distributions. The operators $S_{\hbar}$, $\D_jS_{\hbar}$ and $\Id-2S_{\hbar}$ are bounded self-adjoint operators on $L^2$, with
		\begin{align*}
			\norm{\Id-2S_{\hbar}}_{L^2\to L^2}\leq1,\qquad \norm{\D_m\D_jS_{\hbar}}_{L^2\to L^2}\leq1.
		\end{align*}
		For any $s\in\R$ and $u\in H^s$, we also have
		\begin{align*}
			\norm{S_{\hbar}u}_{H^{s+2}}\leq \hbar^{-2}\norm u_{H^s}.
		\end{align*}
	\end{lemma}
	\begin{proof}
		Multiplication by a real-valued function is self-adjoint on its maximal domain \cite[Section~2.2]{Teschl}. Unitary equivalence under the Fourier transform gives the asserted self-adjointness. The domain identity $\Dom\innerproduct{P_{0,\hbar}}=H^2$ follows from \cite[Section~7.3, Theorem~7.17]{Teschl}. By the Plancherel theorem,
		\begin{align*}
			\norm{\innerproduct{\Id-2S_{\hbar}}u}_{L^2}^2 &=\int_{\R^d}\frac{\innerproduct{\hbar^2\snorm{\vec\xi}^2-1}^2}{\innerproduct{1+\hbar^2\snorm{\vec\xi}^2}^2}\snorm{\widehat u}^2\dd\vec\xi \leq\norm u_{L^2}^2,\\
			\norm{\D_m\D_jS_{\hbar}u}_{L^2}^2 &=\int_{\R^d}\frac{\hbar^4\xi_m^2\xi_j^2}{\innerproduct{1+\hbar^2\snorm{\vec\xi}^2}^2}\snorm{\widehat u}^2\dd\vec\xi \leq\norm u_{L^2}^2.
		\end{align*}
		For $u\in H^s$, we also have
		\begin{align*}
			\norm{S_{\hbar}u}_{H^{s+2}}^2 &=\int_{\R^d}\frac{\innerproduct{1+\snorm{\vec\xi}^2}^2}{\innerproduct{1+\hbar^2\snorm{\vec\xi}^2}^2} \innerproduct{1+\snorm{\vec\xi}^2}^{s}\snorm{\widehat u}^2\dd\vec\xi \leq \hbar^{-4}\norm u_{H^s}^2.
		\end{align*}
		This completes the proof.
	\end{proof}
	
	\begin{lemma}\label{lm:COV}
		Fix $\hbar>0$. For $t>0$ and $\vec z\in\R^d$, define
		\begin{align*}
			G_{t,\hbar}\innerproduct{\vec z} =\innerproduct{4\pi t\hbar^2}^{-\frac{d}{2}}e^{-\frac{\snorm{\vec z}^2}{4t\hbar^2}}.
		\end{align*}
		For $\vec z\ne0$ and $m=1,\ldots,d$, define
		\begin{align*}
			G_\hbar\innerproduct{\vec z} =\int_0^\infty e^{-t}G_{t,\hbar}\innerproduct{\vec z}\dd t,\qquad \mathcal G_{m,\hbar}\innerproduct{\vec z} =-\imagunit \hbar\int_0^\infty e^{-t}\partial_{z_m}G_{t,\hbar}\innerproduct{\vec z}\dd t.
		\end{align*}
		Then $G_\hbar\in W^{1,1}(\R^d)$, $G_\hbar\geq0$, $\norm{G_\hbar}_{L^1}=1$, and $\mathcal G_{m,\hbar}=-\imagunit \hbar\partial_mG_\hbar\in L^1$. For $u\in L^2(\R^d)$, the following identities hold in $L^2$:
		\begin{align}
			\innerproduct{S_{\hbar}u}\innerproduct{\vec x} &=\innerproduct{G_\hbar*u}\innerproduct{\vec x} =\int_{\R^d}G_\hbar\innerproduct{\vec x-\vec y}u\innerproduct{\vec y}\dd\vec y, \label{eq:Kh-convolution}\\
			\innerproduct{\D_mS_{\hbar}u}\innerproduct{\vec x} &=\innerproduct{\mathcal G_{m,\hbar}*u}\innerproduct{\vec x} =\int_{\R^d}\mathcal G_{m,\hbar}\innerproduct{\vec x-\vec y}u\innerproduct{\vec y}\dd\vec y. \label{eq:DKh}
		\end{align}
		Moreover, there is a constant $C_d>0$ independent of $\hbar$ such that
		\begin{align*}
			\norm{\mathcal G_{m,\hbar}}_{L^1}\leq C_d,\qquad \int_{\R^d}\snorm{\vec z}^{\frac{1}{2}}\innerproduct{G_\hbar\innerproduct{\vec z} +\sum_{m=1}^d\snorm{\mathcal G_{m,\hbar}\innerproduct{\vec z}}}\dd\vec z\leq C_d\hbar^{\frac{1}{2}}.
		\end{align*}
	\end{lemma}
	\begin{proof}
		The heat semigroup $e^{t\hbar^2\Delta}$ has Fourier multiplier $e^{-t\hbar^2\snorm{\vec\xi}^2}$. For $u\in L^2$, the Plancherel theorem gives
		\begin{align*}
			\int_0^\infty\norm{e^{-t}e^{t\hbar^2\Delta}u}_{L^2}\dd t &\leq\norm u_{L^2}\int_0^\infty e^{-t}\dd t =\norm u_{L^2}.
		\end{align*}
		The integral therefore converges in $L^2$, and the Fourier transform commutes with this integral. Hence,
		\begin{align*}
			\widehat{\int_0^\infty e^{-t}e^{t\hbar^2\Delta}u\dd t}\innerproduct{\vec\xi} &=\int_0^\infty e^{-t}\widehat{e^{t\hbar^2\Delta}u}\innerproduct{\vec\xi}\dd t\\
			&=\innerproduct{\int_0^\infty e^{-t\innerproduct{1+\hbar^2\snorm{\vec\xi}^2}}\dd t}\widehat u\innerproduct{\vec\xi}\\
			&=\frac{\widehat u\innerproduct{\vec\xi}}{1+\hbar^2\snorm{\vec\xi}^2} =\widehat{S_{\hbar}u}\innerproduct{\vec\xi}.
		\end{align*}
		It follows that
		\begin{align}
			S_{\hbar}u=\int_0^\infty e^{-t}e^{t\hbar^2\Delta}u\dd t\quad\text{in }L^2. \label{eq:heat-representation}
		\end{align}
		Thus, $G_\hbar$ is a rescaled Bessel potential kernel of order two; see \cite[Chapter~II, \S4]{AronszajnSmith}.
		
		The Gaussian heat kernel \cite[\S2.3.1]{EvansPDE} gives, for fixed $t>0$ and almost every $\vec x\in\R^d$,
		\begin{align*}
			\innerproduct{e^{t\hbar^2\Delta}u}\innerproduct{\vec x} =\innerproduct{G_{t,\hbar}*u}\innerproduct{\vec x} =\innerproduct{4\pi t\hbar^2}^{-\frac{d}{2}} \int_{\R^d}e^{-\frac{\snorm{\vec x-\vec y}^2}{4t\hbar^2}}u\innerproduct{\vec y}\dd\vec y.
		\end{align*}
		The one-dimensional Gaussian integral \cite[\S4.3.1]{EvansPDE} yields $\widehat{G_{t,\hbar}}\innerproduct{\vec\xi}=\innerproduct{2\pi}^{-\frac{d}{2}}e^{-t\hbar^2\snorm{\vec\xi}^2}$ and $\int_{\R^d}G_{t,\hbar}\innerproduct{\vec z}\dd\vec z=1$. The value assigned at the origin does not affect the integrals. Since $G_{t,\hbar}\geq0$, Tonelli's theorem \cite[Theorem~2.37]{Folland} gives
		\begin{align*}
			\norm{G_\hbar}_{L^1} &=\int_{\R^d\setminus\setprod0}\int_0^\infty e^{-t}G_{t,\hbar}\innerproduct{\vec z}\dd t\dd\vec z=\int_0^\infty e^{-t}\int_{\R^d\setminus\setprod0}G_{t,\hbar}\innerproduct{\vec z}\dd\vec z\dd t\\
			&=\int_0^\infty e^{-t}\int_{\R^d}G_{t,\hbar}\innerproduct{\vec z}\dd\vec z\dd t =\int_0^\infty e^{-t}\dd t=1.
		\end{align*}
		In particular, $G_\hbar$ is integrable, including near the origin. For $u\in L^2$, Tonelli's theorem and Young's convolution inequality give
		\begin{align*}
			\norm{\int_0^\infty e^{-t}\int_{\R^d} G_{t,\hbar}\innerproduct{\,\cdot-\vec y}\snorm{u\innerproduct{\vec y}}\dd\vec y\dd t}_{L^2}=\norm{G_\hbar*\snorm u}_{L^2} \leq\norm{G_\hbar}_{L^1}\norm u_{L^2}=\norm u_{L^2}.
		\end{align*}
		The nonnegative double integral is therefore finite almost everywhere. Fubini's theorem \cite[Theorem~2.37]{Folland} and \eqref{eq:heat-representation} imply that
		\begin{align}
			S_{\hbar}u\innerproduct{\vec x} &=\int_0^\infty e^{-t}\int_{\R^d} G_{t,\hbar}\innerproduct{\vec x-\vec y}u\innerproduct{\vec y}\dd\vec y\dd t\notag\\
			&=\int_{\R^d}\innerproduct{\int_0^\infty e^{-t}G_{t,\hbar}\innerproduct{\vec x-\vec y}\dd t} u\innerproduct{\vec y}\dd\vec y\notag\\
			&=\int_{\R^d}G_\hbar\innerproduct{\vec x-\vec y}u\innerproduct{\vec y}\dd\vec y =\innerproduct{G_\hbar*u}\innerproduct{\vec x}. \label{eq:Kh-convolution-app}
		\end{align}
		Young's inequality gives $G_\hbar*u\in L^2$, so the identity above holds in $L^2$, i.e., \eqref{eq:Kh-convolution}.
		
		For $a\geq0$, use polar coordinates $\vec z=r\vec\omega$, where $r>0$ and $\vec\omega\in\mathbb S^{d-1}$. The substitution $s=\frac{r^2}{4t\hbar^2}$ gives $\dd r=\hbar\sqrt t\,s^{-\frac{1}{2}}\dd s$, and hence
		\begin{align*}
			\int_{\R^d}\snorm{\vec z}^aG_{t,\hbar}\innerproduct{\vec z}\dd\vec z &=\innerproduct{4\pi t\hbar^2}^{-\frac{d}{2}} \int_{\R^d}\snorm{\vec z}^ae^{-\frac{\snorm{\vec z}^2}{4t\hbar^2}}\dd\vec z\\
			&=\innerproduct{4\pi t\hbar^2}^{-\frac{d}{2}} \int_{\mathbb S^{d-1}}\int_0^\infty r^ae^{-\frac{r^2}{4t\hbar^2}}r^{d-1}\dd r\dd S\innerproduct{\vec\omega}\\
			&=\innerproduct{4\pi t\hbar^2}^{-\frac{d}{2}} \frac{2\pi^{\frac{d}{2}}}{\Gamma\innerproduct{\frac{d}{2}}} \int_0^\infty r^{d+a-1}e^{-\frac{r^2}{4t\hbar^2}}\dd r\\
			&=\innerproduct{4\pi t\hbar^2}^{-\frac{d}{2}} \frac{2\pi^{\frac{d}{2}}}{\Gamma\innerproduct{\frac{d}{2}}} \int_0^\infty\innerproduct{2\hbar\sqrt{ts}}^{d+a-1}e^{-s} \hbar\sqrt t\,s^{-\frac{1}{2}}\dd s\\
			&=\innerproduct{4\pi t\hbar^2}^{-\frac{d}{2}} \frac{2\pi^{\frac{d}{2}}}{\Gamma\innerproduct{\frac{d}{2}}} 2^{d+a-1}\innerproduct{t\hbar^2}^{\frac{d+a}{2}} \int_0^\infty s^{\frac{d+a}{2}-1}e^{-s}\dd s\\
			&=2^a\hbar^at^{\frac{a}{2}}\frac{\Gamma\innerproduct{\frac{d+a}{2}}}{\Gamma\innerproduct{\frac{d}{2}}}.
		\end{align*}
		Tonelli's theorem gives
		\begin{align}
			\int_{\R^d}\snorm{\vec z}^aG_\hbar\innerproduct{\vec z}\dd\vec z &=\int_{\R^d}\snorm{\vec z}^a \innerproduct{\int_0^\infty e^{-t}G_{t,\hbar}\innerproduct{\vec z}\dd t}\dd\vec z\notag\\
			&=\int_0^\infty\int_{\R^d} e^{-t}\snorm{\vec z}^aG_{t,\hbar}\innerproduct{\vec z}\dd\vec z\dd t\notag\\
			&=\int_0^\infty e^{-t}\innerproduct{\int_{\R^d} \snorm{\vec z}^aG_{t,\hbar}\innerproduct{\vec z}\dd\vec z}\dd t\notag\\
			&=2^a\hbar^a\frac{\Gamma\innerproduct{\frac{d+a}{2}}}{\Gamma\innerproduct{\frac{d}{2}}} \int_0^\infty e^{-t}t^{\frac{a}{2}}\dd t\notag\\
			&=2^a\hbar^a \frac{\Gamma\innerproduct{1+\frac{a}{2}}\Gamma\innerproduct{\frac{d+a}{2}}} {\Gamma\innerproduct{\frac{d}{2}}}. \label{eq:moment-app}
		\end{align}
		Taking $a=\frac{1}{2}$ yields $\int_{\R^d}\snorm{\vec z}^{\frac{1}{2}}G_\hbar\innerproduct{\vec z}\dd\vec z=C_d\hbar^{\frac{1}{2}}$.
		
		To estimate the derivative kernel, we differentiate $G_{t,\hbar}$ directly:
		\begin{align*}
			\partial_{z_m}G_{t,\hbar}\innerproduct{\vec z} =-\frac{z_m}{2t\hbar^2}\innerproduct{4\pi t\hbar^2}^{-\frac{d}{2}} e^{-\frac{\snorm{\vec z}^2}{4t\hbar^2}} =-\frac{z_m}{2t\hbar^2}G_{t,\hbar}\innerproduct{\vec z}.
		\end{align*}
		For $a\geq0$, the preceding moment formula gives
		\begin{align*}
			\hbar\int_{\R^d}\snorm{\vec z}^a\snorm{\partial_{z_m}G_{t,\hbar}\innerproduct{\vec z}}\dd\vec z &=\frac1{2t\hbar}\int_{\R^d}\snorm{\vec z}^a\snorm{z_m}G_{t,\hbar}\innerproduct{\vec z}\dd\vec z\\
			&\leq\frac1{2t\hbar}\int_{\R^d}\snorm{\vec z}^{a+1}G_{t,\hbar}\innerproduct{\vec z}\dd\vec z =2^a\hbar^at^{\frac{a-1}{2}} \frac{\Gamma\innerproduct{\frac{d+a+1}{2}}}{\Gamma\innerproduct{\frac{d}{2}}}.
		\end{align*}
		In particular,
		\begin{align*}
			\hbar\int_0^\infty e^{-t}\norm{\partial_{z_m}G_{t,\hbar}}_{L^1}\dd t &\leq\frac{\Gamma\innerproduct{\frac{d+1}{2}}}{\Gamma\innerproduct{\frac{d}{2}}} \int_0^\infty e^{-t}t^{-\frac{1}{2}}\dd t =\frac{\sqrt\pi\,\Gamma\innerproduct{\frac{d+1}{2}}}{\Gamma\innerproduct{\frac{d}{2}}}<\infty.
		\end{align*}
		Thus, $\int_0^\infty e^{-t}\partial_{z_m}G_{t,\hbar}\dd t$ converges in $L^1$. For $\varphi\in C_c^\infty$, Fubini's theorem and integration by parts give
		\begin{align*}
			\int_{\R^d}G_\hbar\innerproduct{\vec z}\partial_{z_m}\varphi\innerproduct{\vec z}\dd\vec z &=\int_0^\infty e^{-t}\int_{\R^d}G_{t,\hbar}\innerproduct{\vec z}\partial_{z_m}\varphi\innerproduct{\vec z}\dd\vec z\dd t\\
			&=-\int_0^\infty e^{-t}\int_{\R^d}\partial_{z_m}G_{t,\hbar}\innerproduct{\vec z}\varphi\innerproduct{\vec z}\dd\vec z\dd t\\
			&=-\int_{\R^d}\innerproduct{\int_0^\infty e^{-t}\partial_{z_m}G_{t,\hbar}\innerproduct{\vec z}\dd t}\varphi\innerproduct{\vec z}\dd\vec z.
		\end{align*}
		It follows that $G_\hbar\in W^{1,1}$ and
		\begin{align*}
			\mathcal G_{m,\hbar}\innerproduct{\vec z} =-\imagunit \hbar\partial_{z_m}G_\hbar\innerproduct{\vec z} =-\imagunit \hbar\int_0^\infty e^{-t}\partial_{z_m}G_{t,\hbar}\innerproduct{\vec z}\dd t\quad\text{in }L^1.
		\end{align*}
		For any $u\in L^2$, differentiation of the convolution in the weak sense gives
		\begin{align*}
			\innerproduct{\D_mS_{\hbar}u}\innerproduct{\vec x} &=-\imagunit \hbar\partial_{x_m}\int_{\R^d}G_\hbar\innerproduct{\vec x-\vec y}u\innerproduct{\vec y}\dd\vec y\\
			&=\int_{\R^d}\mathcal G_{m,\hbar}\innerproduct{\vec x-\vec y}u\innerproduct{\vec y}\dd\vec y =\innerproduct{\mathcal G_{m,\hbar}*u}\innerproduct{\vec x}.
		\end{align*}
		The right-hand side belongs to $L^2$ by Young's inequality. Thus, \eqref{eq:DKh} holds in $L^2$. Moreover,
		\begin{align}
			\int_{\R^d}\snorm{\vec z}^a\snorm{\mathcal G_{m,\hbar}\innerproduct{\vec z}}\dd\vec z &\leq \hbar\int_0^\infty e^{-t}\int_{\R^d} \snorm{\vec z}^a\snorm{\partial_{z_m}G_{t,\hbar}\innerproduct{\vec z}}\dd\vec z\dd t\notag\\
			&\leq2^a\hbar^a\frac{\Gamma\innerproduct{\frac{d+a+1}{2}}}{\Gamma\innerproduct{\frac{d}{2}}} \int_0^\infty e^{-t}t^{\frac{a-1}{2}}\dd t\notag\\
			&=2^a\hbar^a \frac{\Gamma\innerproduct{\frac{a+1}{2}}\Gamma\innerproduct{\frac{d+a+1}{2}}} {\Gamma\innerproduct{\frac{d}{2}}},\qquad a\geq0. \label{eq:derivative-moment-app}
		\end{align}
		Taking $a=0$ and $a=\frac{1}{2}$, respectively, gives
		\begin{align*}
			\norm{\mathcal G_{m,\hbar}}_{L^1}\leq C_d,\qquad \int_{\R^d}\snorm{\vec z}^{\frac{1}{2}}\snorm{\mathcal G_{m,\hbar}\innerproduct{\vec z}}\dd\vec z \leq C_d\hbar^{\frac{1}{2}}.
		\end{align*}
		Combining these bounds with \eqref{eq:moment-app} for $a=\frac{1}{2}$ and summing over $m$ completes the proof.
	\end{proof}
	
	\section{The stretched Green kernel and interior decay}
	\label{app:green-truncation}
	\begin{lemma}
		\label{lem:green-kernel}
		For $\vec x\ne\vec y$, define
		\begin{align*}
			s_j=\int_{y_j}^{x_j}\sigma\innerproduct t\dd t, \qquad \widetilde x_j-\widetilde y_j=x_j-y_j+\zeta s_j+\imagunit s_j.
		\end{align*}
		Then
		\begin{align}
			\Ima\rho\innerproduct{\widetilde{\vec x},\widetilde{\vec y}}&\geq\frac1{\sigma_0\snorm{\vec x-\vec y}}\sum_{j=1}^d\snorm{s_j}^2, \label{eq:rho-imaginary}\\
			c_d\snorm{\vec x-\vec y}&\leq\snorm{\rho\innerproduct{\widetilde{\vec x},\widetilde{\vec y}}}\leq M\snorm{\vec x-\vec y}, \qquad \text{where}~~c_d=\begin{cases} \innerproduct{1+\sigma_0}^{-1}, & d=2, \\
				\sqrt{\cos\omega_0}, & d=3. \end{cases} \label{eq:rho-comparable}
		\end{align}
		The kernel $G_d$ satisfies the fundamental-solution identity
		\begin{align}
			\midproduct{-\grad_{\vec x}\cdot\innerproduct{\mat A\innerproduct{\vec x}\grad_{\vec x}}-k^2J\innerproduct{\vec x}}G_d\innerproduct{\widetilde{\vec x},\widetilde{\vec y}}=\delta_{\vec y}. \label{eq:green-fundamental}
		\end{align}
		Moreover, there exists $C_G>0$ independent of $k$, $\vec x$ and $\vec y$ such that, for $j=0,1$,
		\begin{align}
			\snorm{\grad_{\vec x}^{j}G_2\innerproduct{\widetilde{\vec x},\widetilde{\vec y}}}&\leq C_Gk^j\midproduct{\innerproduct{k\snorm{\vec x-\vec y}}^{-\frac{1}{2}}+j\innerproduct{k\snorm{\vec x-\vec y}}^{-j}}e^{-\frac{k\Ima\rho}{2}}, \label{eq:G2-pointwise}\\
			\snorm{\grad_{\vec x}^{j}G_3\innerproduct{\widetilde{\vec x},\widetilde{\vec y}}}&\leq C_G\midproduct{k^j\snorm{\vec x-\vec y}^{-1}+j\snorm{\vec x-\vec y}^{-j-1}}e^{-k\Ima\rho}. \label{eq:G3-pointwise}
		\end{align}
	\end{lemma}
	\begin{proof}
		Here \eqref{eq:rho-imaginary} is proved in \cite[Lemma~2.2]{BaiDuWu}. Its proof uses only $\zeta\geq0$ and $0\leq\sigma\leq\sigma_0$, and requires no upper bound on $\zeta$. For $d=2$, the corresponding arguments for \eqref{eq:rho-comparable}, \eqref{eq:G2-pointwise}, and the fundamental-solution identity are given in \cite[Lemmas~2.1--2.3]{WangZheng}.
		
		For $d=3$, the following argument applies to any fixed $\zeta\geq\frac{1}{\sqrt3}$. Set $r=\snorm{\vec x-\vec y}$ and $t_j=\frac{\snorm{x_j-y_j}^2}{r^2}$. Let $m_j=\frac{s_j}{x_j-y_j}$ when $x_j\ne y_j$, and $m_j=0$ otherwise. Then $0\leq m_j\leq\sigma_0$, $\sum_jt_j=1$, and
		\begin{align*}
			\rho^2=r^2\sum_{j=1}^3t_j\midproduct{1+\innerproduct{\zeta+\imagunit}m_j}^2.
		\end{align*}
		From Lemma~\ref{lem:sectorial-A}, writing $1+\innerproduct{\zeta+\imagunit}m_j=r_je^{\imagunit\theta_j}$, we have $0\leq\theta_j\leq\omega_0$, and hence
		\begin{align*}
			\snorm{\sum_{j=1}^3t_jr_j^2e^{2\imagunit\theta_j}}\geq\Rea\innerproduct{e^{-\imagunit\omega_0}\sum_{j=1}^3t_jr_j^2e^{2\imagunit\theta_j}}=\sum_{j=1}^3t_jr_j^2\cos\innerproduct{2\theta_j-\omega_0}\geq\cos\omega_0.
		\end{align*}
		This proves the lower bound in \eqref{eq:rho-comparable} for $d=3$. The upper bound follows from $\snorm{1+\innerproduct{\zeta+\imagunit}m_j}\leq M$. Differentiating $\frac{e^{\imagunit k\rho}}{4\pi\rho}$ and using \eqref{eq:rho-comparable} gives \eqref{eq:G3-pointwise}; see \cite[Lemmas~2.2--2.4]{BaiDuWu}.
	\end{proof}
	\begin{lemma}
		\label{lem:green-interior-decay}
		Let $z\in H^1(B_{L+D_c})$ satisfy \eqref{eq:green-z-homogeneous}. If $D_c\geq\bar l$ and $k\geq k_0$, then
		\begin{align}
			\enorm z_{B_L}\leq\begin{cases} C_GM\innerproduct{kD_c}^{\frac{1}{2}}e^{-\frac{k\sigma_0D_c}{40\sqrt2}}\enorm z_{B_{L+D_c}}, & d=2, \\[2mm]
				C_GM^2\innerproduct{kD_c}e^{-\frac{k\sigma_0D_c}{20\sqrt3}}\enorm z_{B_{L+D_c}}, & d=3. \end{cases} \label{eq:green-interior-decay}
		\end{align}
	\end{lemma}
	\begin{proof}
		Increasing $k_0$ if necessary, we assume that $k_0\bar l\geq2$, so $kD_c\geq2$. Choose a product cutoff $\chi\innerproduct{\vec x}=\prod_{j=1}^d\chi_j\innerproduct{x_j}\in C_c^\infty(B_{L+D_c})$ such that $0\leq\chi_j\leq1$, $\chi_j=1$ for $\snorm{x_j}\leq L+\frac{D_c}{2}$, and $\chi_j=0$ for $\snorm{x_j}\geq L+D_c$. Then
		\begin{align*}
			\snorm{\partial_j\chi}\leq CD_c^{-1}, \qquad \snorm{\partial_{jj}\chi}\leq CD_c^{-2}, \qquad j=1,\ldots,d,
		\end{align*}
		and $\partial_jA_{jj}=0$ on $\operatorname{supp}\partial_j\chi$. Consequently, $\innerproduct{\partial_jA_{jj}}\partial_j\chi=0$ for $j=1,\ldots,d$.
		
		Extending $\chi z$ by zero to $\R^d$, we obtain from \eqref{eq:green-z-homogeneous}
		\begin{align*}
			\innerproduct{-\divg\innerproduct{\mat A\grad}-k^2J}\innerproduct{\chi z}&=-z\divg\innerproduct{\mat A\grad\chi}-2\mat A\grad\chi\cdot\grad z=:-g_\chi.
		\end{align*}
		Using $\innerproduct{\partial_jA_{jj}}\partial_j\chi=0$, the coefficient bounds \eqref{eq:truncated-coefficient-bounds}, and $kD_c\geq2$, we obtain
		\begin{align}
			\norm{g_\chi}_{L^2}
			&\leq\norm{z\divg\innerproduct{\mat A\grad\chi}}_{L^2}+2\norm{\mat A\grad\chi\cdot\grad z}_{L^2}\notag\\
			&\leq CM^{d-1}\innerproduct{D_c^{-2}\norm z_{L^2(B_{L+D_c})}+D_c^{-1}\norm{\grad z}_{L^2(B_{L+D_c})}}\notag\\
			&\leq CM^{d-1}D_c^{-1}\enorm z_{B_{L+D_c}}. \label{eq:gchi-bound}
		\end{align}
		The Green representation gives, for almost every $\vec x\in B_L$,
		\begin{align}
			z\innerproduct{\vec x}=\chi\innerproduct{\vec x}z\innerproduct{\vec x}=-\int_{\R^d}G_d\innerproduct{\widetilde{\vec x},\widetilde{\vec y}}g_\chi\innerproduct{\vec y}\dd\vec y. \label{eq:green-representation}
		\end{align}
		If $\vec y\in\operatorname{supp}g_\chi$, then $\snorm{y_j}\geq L+\frac{D_c}{2}$ for at least one index $j$. Since $\sigma=\sigma_0$ for $\snorm t\geq\bar l$, any $\vec x\in B_L$ satisfies
		\begin{align*}
			\snorm{\vec x-\vec y}\geq\frac{D_c}{2}, \qquad \snorm{\vec x-\vec y}\leq\sqrt d\innerproduct{2L+D_c}, \qquad \snorm{\int_{y_j}^{x_j}\sigma\innerproduct t\dd t}\geq\frac{\sigma_0D_c}{2}.
		\end{align*}
		Since $D_c\geq\bar l$, we have $L=\bar l+D_c\leq2D_c$. Estimate \eqref{eq:rho-imaginary} therefore yields
		\begin{align}
			\Ima\rho\innerproduct{\widetilde{\vec x},\widetilde{\vec y}}\geq\frac{\sigma_0D_c^2}{4\sqrt d\innerproduct{2L+D_c}}\geq\frac{\sigma_0D_c}{20\sqrt d}. \label{eq:rho-separated}
		\end{align}
		The kernel bounds \eqref{eq:G2-pointwise}--\eqref{eq:G3-pointwise} give, for $j=0,1$,
		\begin{align}
			\innerproduct{\int_{B_L}\int_{\operatorname{supp}g_\chi}\snorm{\grad_{\vec x}^{j}G_2\innerproduct{\widetilde{\vec x},\widetilde{\vec y}}}^2\dd\vec y\dd\vec x}^{\frac{1}{2}}&\leq C_GD_c^{\frac{3}{2}}k^{j-\frac{1}{2}}e^{-\frac{k\sigma_0D_c}{40\sqrt2}}, \label{eq:G2-HS}\\
			\innerproduct{\int_{B_L}\int_{\operatorname{supp}g_\chi}\snorm{\grad_{\vec x}^{j}G_3\innerproduct{\widetilde{\vec x},\widetilde{\vec y}}}^2\dd\vec y\dd\vec x}^{\frac{1}{2}}&\leq C_GD_c^2k^je^{-\frac{k\sigma_0D_c}{20\sqrt3}}. \label{eq:G3-HS}
		\end{align}
		
		Applying Cauchy--Schwarz to \eqref{eq:green-representation} and using \eqref{eq:G2-HS}--\eqref{eq:G3-HS} and \eqref{eq:gchi-bound}, we obtain
		\begin{align*}
			\enorm z_{B_L}&\leq C_GM\innerproduct{kD_c}^{\frac{1}{2}}e^{-\frac{k\sigma_0D_c}{40\sqrt2}}\enorm z_{B_{L+D_c}}, && d=2,\\
			\enorm z_{B_L}&\leq C_GM^2\innerproduct{kD_c}e^{-\frac{k\sigma_0D_c}{20\sqrt3}}\enorm z_{B_{L+D_c}}, && d=3.
		\end{align*}
		This proves \eqref{eq:green-interior-decay}. For fixed $\sigma_0>0$ and $M$, as $t\to\infty$,
		\begin{align*}
			t^{\frac{1}{2}}e^{-\frac{\sigma_0t}{40\sqrt2}}\longrightarrow0, \qquad te^{-\frac{\sigma_0t}{20\sqrt3}}\longrightarrow0.
		\end{align*}
		Hence, the factor on the right-hand side of \eqref{eq:green-interior-decay} is at most $2^{-\frac{d}{2}-1}$ for sufficiently large $k\sigma_0D_c$.
	\end{proof}
	
	\bibliographystyle{siam}
	\bibliography{UPML}
\end{document}